\documentclass[twocolumn]{article}
\usepackage[utf8]{inputenc}
\usepackage{mathrsfs}
\usepackage{amsmath,amssymb,amstext,amsfonts, amsthm}
\usepackage{mathtools,cuted}
\usepackage{algorithm}
\usepackage{algpseudocode}
\usepackage{geometry}
\usepackage{multirow}
\usepackage{tikz}
\usepackage{xcolor}
\usetikzlibrary{quantikz}
\usepackage{sectsty}
\subsectionfont{\normalfont\itshape \large}
\subsubsectionfont{\normalfont\itshape \large}
\sectionfont{\large}
\usepackage{subcaption}
\usepackage{pgfplots}
\pgfplotsset{compat = newest}

\providecommand{\keywords}[1]
{
	
  \textbf{Keywords} : #1
}
\newcommand{\C}{\mathbb{C}}
\newcommand{\R}{\mathbb{R}}

\newcommand{\qH}{\mathbb{H}_q}

\newcommand{\mal}{\otimes}

\newcommand{\mlqae}{\text{ MLQAE }}
\newcommand{\return}{\State \textbf{Return }}

\title{QFT-based Homogenization}
\author{F. Givois$^1$ \and M. Kabel$^1$ \and N. Gauger$^2$ }
\date{$^1$\emph{Fraunhofer ITWM, Kaiserslautern, Germany}\\
    $^2$\emph{Chair for scientific computing, University of Kaiserslautern, Germany}\\[2ex]%
    \today
}

\newcommand{\opbf}{data/simulator/1D/phase/bell_state/qft/qb4.dat}
\newcommand{\opbt}{data/simulator/1D/phase/bell_state/qft/qb2.dat}
\newcommand{\opbtr}{data/simulator/1D/phase/bell_state/qft/qb3.dat}
\newcommand{\opbfi}{data/simulator/1D/phase/bell_state/qft/qb5.dat}

\newcommand{\opbLt}{data/simulator/1D/phase/bell_state/L2/L22.dat}

\newcommand{\Yopbt}{data/real_devices/ibmq_ehningen/phase/bell_state/qft/qb2.dat}
\newcommand{\Yopbth}{data/real_devices/ibmq_ehningen/phase/bell_state/qft/qb3.dat}
\newcommand{\YoHbt}{data/real_devices/ibmq_ehningen/Hadamard/bell_state/qft/qb2.dat}
\newcommand{\YopbLt}{data/real_devices/ibmq_ehningen/phase/bell_state/L2/L22.dat}
\newcommand{\YopbLth}{data/real_devices/ibmq_ehningen/phase/bell_state/L2/L23.dat}

\begin{document}

\twocolumn[
  \begin{@twocolumnfalse}
    \maketitle
    \begin{abstract}
      Efficient numerical characterization is a key problem in composite material analysis. To follow accuracy improvement in image tomography, memory efficient methods of numerical characterization have been developed. Among them, an FFT based solver has been proposed by Moulinec and Suquet (1994,1998) bringing down numerical characterization complexity to the FFT complexity. Nevertheless, recent development of tomography sensors made memory requirement and calculation time reached another level. To avoid this bottleneck, the new leaps in the field of Quantum Computing have been used. This paper will present the application of the Quantum Fourier Transform (QFT) to replace the Fast Fourier Transform (FFT) in Moulinec and Suquet algorithm. It will mainly focused on how to read out Fourier coefficients stored in a quantum state. First, a reworked Hadamard test algorithm applied with Most likelihood amplitude estimation (MLQAE) is used to determine the quantum coefficients. Second, an improvement avoiding Hadamard test is presented in case of Material characterization on mirrored domain. Finally, this last algorithm is applied to Material homogenization to determine effective stiffness of basic geometries. 
    \end{abstract} \vspace{10pt}
    \vspace{20pt}
  \end{@twocolumnfalse}
]

\section{Introduction}

The computation of effective material laws for composites with complicated micro-structures is a difficult task. As analytical description of composite materials exists only for specific classes of materials (short fiber reinforced plastic for example), the best method is to approximate material effective behaviour by simulating material unit cell on the microscopic scale, based on a tomography image. To do so, the approximation of the solution of an homogenization problem described by an elliptic PDE is needed. As the computation domain can be very wide, a memory efficient FFT-based homogenization algorithm has been developed by Moulinec and Suquet \cite{MOULINEC199869}. This algorithm solves the PDE iteratively with a matrix-free gradient-descent method. But despite the fact that this algorithm is matrix free, it can be very costly as it involves 3D Fourier transforms of large data domain at each iteration to invert the preconditioner. During the last decades, the computation power needed for material characterization skyrockets with the improvement of material tomography imaging. This improvement in the accuracy of material images turned the computation domains to be very wide (more than a terabyte of data), this leads to way longer computation time due to Fourier transform complexity, and to memory bottleneck. 
\\

On the other hand, introduced by Richard Feynman in the 80's \cite{feynman2018simulating} Quantum computing wildly developed during the past 40 years in divert fields of research. Recent leaps in the development of Quantum computers initially introduce by companies such as IBM or Google seem to let the door open for practical application of quantum algorithms. Especially, development of programming interfaces and quantum calculation platforms such as IBM Qiskit library and NISQ resulted in the development of various applied quantum algorithms with in the long term application to real world problems. Variational Quantum algorithms (VQE) \cite{Peruzzo_2014} \cite{Kandala_2017} are for example having applications in fields such as Quantum Chemistry or Quantum Machine Learning  
\\

In this paper we study an implementation of an adaptation of Moulinec and Suquet algorithm \cite{MOULINEC199869} for symmetric inputs developed by M.K. and H.Grimm-Strele \cite{kabel2016mixed} using a replacement of the Fast Fourier Transform (FFT) by a Quantum Fourier Transform (QFT). The quantum coefficients are estimated with an algorithm mixing Hadamard test \cite{aharonov2006polynomial} with most likelyhood quantum amplitude estimation algorithm (MLQAE) \cite{suzuki2020amplitude}. As a first intermediate result we focused on MLQAE method in order to determine a way to estimate \emph{a priori} the number of amplification necessary for a desired accuracy. Secondly, we describe a more efficient quantum coefficients estimation method for symmetric boundary conditions which avoid Hadamard test algorithm and its limitations. As a conclusion we present the overall QFT-based Homogenization algorithm and some results for trivial geometries obtained on IBMQ quantum computer QASM simulator. Despite their apparent complexity inefficiency, our algorithms are used as a proof of concept for the practical use of QFT in hybrid QC-HPC algorithms. In this goal, we validate our algorithm on QASM simulator using different geometries with known analytical solutions.

\section{FFT based Homogenization}

\subsection{Periodic Boundary conditions}

Consider the 2 dimensional square $Y=[0,L_1] \times [0,L_2]\in \mathbb{R}^2$ of edge length $L_i \; (i=1,2)$. Let $\boldsymbol{\sigma} : Y \xrightarrow{} \mathbb{R}^{3}$ the local stress field, $\boldsymbol{u} : Y \xrightarrow{} \mathbb{R}^2$ the local displacement, the local strain field $\boldsymbol{\varepsilon} : \mathbb{R}^2 \xrightarrow{} \mathbb{R}^{3}$ is defined such as $\boldsymbol{\varepsilon}(\boldsymbol{u}(\boldsymbol{x}))= \boldsymbol{E}+\boldsymbol{\varepsilon}(\boldsymbol{u^*}(\boldsymbol{x})) \; \forall \boldsymbol{x} \in Y$ where $\boldsymbol{E}$ is the average strain and $\boldsymbol{\varepsilon}(\boldsymbol{u^*}(\boldsymbol{x}))$ a fluctuation term. The stress behavior of an elastic composite material of local stiffness $\boldsymbol{c}(\boldsymbol{x}) : \R^{2} \xrightarrow{} \mathbb{R}^{3 \times 4}$ can be described by generalized Hooke's law : $\boldsymbol{\sigma}(\boldsymbol{x})= \boldsymbol{C}(\boldsymbol{x}) : \boldsymbol{\varepsilon}(\boldsymbol{u}(\boldsymbol{x})), \; \forall \boldsymbol{x} \in Y$. 
The equation system to solve is deduced from the equation of motion 
\begin{equation}
\label{eqn:system}
    \begin{cases} \boldsymbol{\nabla} \boldsymbol{\sigma}(\boldsymbol{x}) = \boldsymbol{0}, \forall \boldsymbol{x} \in Y, \\ \boldsymbol{\sigma}(\boldsymbol{x})= \boldsymbol{C}(\boldsymbol{x}) : \boldsymbol{\varepsilon}(\boldsymbol{u^*}(\boldsymbol{x}))+\boldsymbol{E}, \; \forall \boldsymbol{x} \in Y. \end{cases}
\end{equation}
To solve this system a method has been  developed by Moulinec and Suquet \cite{MOULINEC199869} in the case of periodic boundary conditions for the displacement.

The algorithm is based on solving the equivalent Lippmann and Schwinger equation for elasticity using a polarization operator $\boldsymbol{C}^0 \in \mathbb{R}^{3 \times 4}$ and a green operator $\boldsymbol{\Gamma}^0 : Y \xrightarrow{} \mathbb{R}^{2 \times 2 \times 4}$ known explicitly in Fourier space. On a sampled grid of $N=N_1 \times N_2=2^{n_1} \times 2^{n_2}$ pixels labeled by $\boldsymbol{h}=\begin{pmatrix} h_1 \\ h_2\end{pmatrix}$, $\boldsymbol{x}(\boldsymbol{h})=\Bigg((h_1-1)\frac{L_1}{N_1},(h_2-1)\frac{L_2}{N_2}\Bigg)$ where $ h_1 = 1,...,N_1,h_2 = 1,...,N_2$.
The algorithm can be written,
\begin{algorithm}[H]
\caption{STRAIN}\label{alg:FFT_based_alg}
\begin{algorithmic}
\Require $\boldsymbol{E},\boldsymbol{C}$
\State $\boldsymbol{\varepsilon}^0(\boldsymbol{x})=\boldsymbol{E}, \forall \boldsymbol{x} \in Y,$
\State $k=1$
\While{convergence is not reached}
\State $\boldsymbol{\tau}^k = (\boldsymbol{C}(\boldsymbol{x})-\boldsymbol{C}^0):\boldsymbol{\varepsilon}^k(\boldsymbol{x}) $
\State $\boldsymbol{\hat{\tau}}^k(\boldsymbol{\xi})=\text{FFT}(\boldsymbol{\tau}^k(\boldsymbol{x})),$

\State $\boldsymbol{\hat{\varepsilon}}^{k+1}(\boldsymbol{\xi})=\boldsymbol{\hat{\Gamma}}^0(\boldsymbol{\xi}):\boldsymbol{\hat{\tau}}^k(\boldsymbol{\xi}),$
\State $\boldsymbol{\varepsilon}^{k+1}(\boldsymbol{x})=\boldsymbol{E}-\text{FFT}^{-1}(\boldsymbol{\hat{\varepsilon}}^{k+1}(\boldsymbol{\xi})),$
\State $k=k+1$,
\State convergence test,
\EndWhile
\State $\boldsymbol{\sigma}(\boldsymbol{x})=\boldsymbol{C}(\boldsymbol{x}):\boldsymbol{\varepsilon}^{k}(\boldsymbol{x}),$
\return $\boldsymbol{\sigma}$
\end{algorithmic}
\end{algorithm}
In the algorithm $\text{FFT}$ describes the Fast Fourier transform. Fourier coefficients of green operator are given by
\begin{equation}
\begin{aligned}
    \boldsymbol{\hat{\Gamma}}^0_{jklm}(\boldsymbol{\xi})= &\frac{1}{4\mu^0 \|\boldsymbol{\xi}\|^2}(\delta_{lj}\xi_m\xi_k+\delta_{mj}\xi_l\xi_k+\delta_{lk}\xi_j\xi_m+\delta_{km}\xi_l\xi_j)\\& -\frac{\lambda^0+\mu^0}{\mu^0(\lambda^0+2\mu^0)}\frac{\xi_j\xi_k\xi_l\xi_m}{\|\boldsymbol{\xi}\|^4}, \quad i,j,k,l = 1,2
    \end{aligned}
\end{equation}
where $\mu^0$ and $\lambda^0$ are the lame coefficients of a reference material chosen according to \cite{https://doi.org/10.1002/nme.275}. The $\delta$ operator is defined such as $\delta_{kl} = \begin{cases} 1 \text{ if } k=l, \\ 0 \text{ otherwise}. \end{cases} $ The $\boldsymbol{\xi}$ are the frequencies in Fourier space, on an $N=N_1 \times N_2=2^{n_1} \times 2^{n_2}$ grid labeled by $\boldsymbol{h}=\begin{pmatrix} h_1 \\ h_2\end{pmatrix}$, $\boldsymbol{\xi}(\boldsymbol{h})=  \Bigg((-\frac{N_1}{2}+h_1)\frac{1}{L_1},(-\frac{N_2}{2}+h_2)\frac{1}{L_2}\Bigg), \quad h_1 = 1,...,N_1,h_2 = 1,...,N_2.$   Discrete error used to test convergence is given by \begin{equation}
    e^{k+1}=\sqrt{\frac{\|\boldsymbol{\varepsilon}^{k+1}
    -\boldsymbol{{\varepsilon}}^k\|^2 }{\|\boldsymbol{E}\|^2}}.
\end{equation}
A step of this algorithm has a complexity bounded by the complexity of the Fast Fourier Transform which is in $\mathcal{O}(N\log(N))$. Nevertheless a naive implementation of the discrete Fourier transform would have given a complexity of $\mathcal{O}(N^2)$, complexity reduction is due to the divide and conquer principle used in the original paper of Cooley and Tukey in 1965 \cite{Cooley1965AnAF}.\\

\subsection{Effective stiffness}

In case of linear elastic materials, the interesting physical quantity to determine is the macroscopic stiffness tensor $\boldsymbol{C^*} \in \mathbb{R}^{3 \times 4}$ which links the average stress field $\langle\boldsymbol{\sigma}\rangle_Y$ and the average strain field $\langle\boldsymbol{\varepsilon}\rangle_Y$ according to the Hill-Mandel condition\begin{equation}
    \langle\boldsymbol{\sigma}\rangle_Y=\boldsymbol{C^*}\langle\boldsymbol{\varepsilon}\rangle_Y.
\end{equation}
Noticing that $\langle\boldsymbol{\varepsilon}\rangle_Y= \boldsymbol{E}$ is prescribed in the FFT-based algorithm, a general homogenization algorithm can be deduced from this relation (using de voigt notations) :
\begin{algorithm}[H]
\caption{QFT-based Homogenization}
\begin{algorithmic}
\Require $\boldsymbol{C}$
\For{$j=1,3$}
\State $\boldsymbol{E}=\boldsymbol{0}$
\State $\boldsymbol{E}_{j}=\boldsymbol{1}$
\State $\boldsymbol{\sigma}= \text{STRAIN}(\boldsymbol{E},\boldsymbol{C})$
\For{$k=1,3$}
\State $C^*_{jk}=\langle\boldsymbol{\sigma}_k\rangle_Y$
\EndFor
\EndFor
\return $\boldsymbol{C^*}$.
\end{algorithmic}
\end{algorithm}
The complexity class of this algorithm is the same of the complexity class of the strain calculation algorithm which is in $\mathcal{O}(N\log(N))$. The question is now to see if it's possible to take advantage of quantum computing methods and more especially Quantum Fourier Transform (QFT). for Periodic mixed boundary conditions the average value is taken on the initial domain and not on the mirrored one.

\subsection{Mixed Uniform Boundary conditions}

To describe effective material behavior, a set of different boundary conditions (BC) can be used. First, Periodic Boundary Conditions (PBC), it has the pro to be very simple to implement using the Moulinec and Suquet algorithm as it is but are not sufficient for non periodic structures according to Chevalier \cite{chevalier2007validation}. Following this idea Periodic-compatible Mixed Uniform Boundary Conditions (PMUBC) have been introduced by Wiegmann \cite{wiegmann1999fast} and Pahr \& Zysset \cite{pahr2008influence}. A remarkable equivalence between PMUBC and Periodic Boundary Conditions on mirrored domain has been used by Bakhalov \cite{bakhvalov2012homogenisation} to reduce computational domain and so improve efficiency for periodic boundary conditions, furthermore, a PMUBC Homogenization solver using PBC on mirrored domain has been developed by Grimm Strele \& Kabel \cite{GrimmStrele2021FFTB}. Using mirrored domain as on Figure \ref{fig:mirrored_domain} to solve Equation \ref{eqn:system} has to be done wisely using different symmetries for the strain fields in function of the initial loading  \cite{GrimmStrele2021FFTB}, it can be summarized in 2 dimensions in Table \ref{table:1} using de Voigt notations. 
\begin{table}[H]
\centering
\begin{tabular}{ |c|c c c| } 
\hline
load case & $\varepsilon_{11}$ & $\varepsilon_{22}$ & $\varepsilon_{12}$ \\
\hline
$xx$ & EE & EE & OO \\
$yy$ & EE & EE  & OO  \\
$xy$ & OO & OO & EE \\
\hline
\end{tabular}
\caption{Table of strain tensor symmetries depending on strain load spatial direction. "E" means even symmetry and "O" means odd symmetry.  }
\label{table:1}
\end{table}
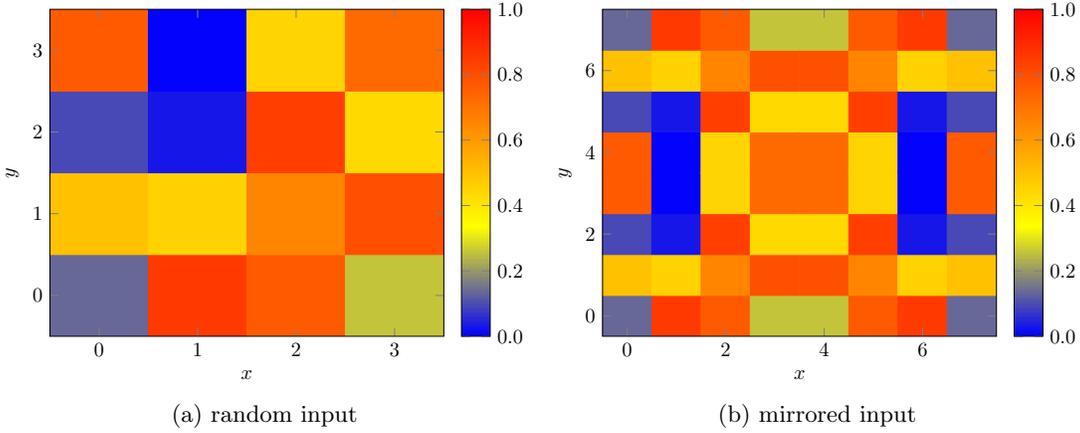
\begin{figure*}
    \centering
\begin{subfigure}{0.4\textwidth}
    \centering
        \resizebox{1\textwidth}{!}{
        \begin{tikzpicture}
          \begin{axis}[
            view={0}{90},   
            xlabel=$x$,
            ylabel=$y$,
            ymin=-0.5,ymax=3.5,
            xmin=-0.5,xmax=3.5,
            colorbar,
            colorbar style={
                yticklabel style={
                    /pgf/number format/.cd,
                    fixed,
                    precision=1,
                    fixed zerofill,
                },
            },
            enlargelimits=false,
            axis on top,
            point meta min=0,
            point meta max=1,
        ]
            \addplot [matrix plot*,point meta=explicit] file {data/simulator/2D/random3/ini_vec.dat};
          \end{axis}
        \end{tikzpicture}
        }
        \caption{random input}
        \label{fig:my_label}
\end{subfigure}
\begin{subfigure}{0.4\textwidth}
    \centering
        \resizebox{1\textwidth}{!}{
        \begin{tikzpicture}
          \begin{axis}[
            view={0}{90},   
            xlabel=$x$,
            ylabel=$y$,
            ymin=-0.5,ymax=7.5,
            xmin=-0.5,xmax=7.5,
            colorbar,
            colorbar style={
                yticklabel style={
                    /pgf/number format/.cd,
                    fixed,
                    precision=1,
                    fixed zerofill,
                },
            },
            enlargelimits=false,
            axis on top,
            point meta min=0,
            point meta max=1,
        ]
            \addplot [matrix plot*,point meta=explicit] file {data/simulator/2D/random3/mirrored_ini_vec.dat};
          \end{axis}
        \end{tikzpicture}
        }
        \caption{mirrored input}
        \label{fig:my_label}
\end{subfigure}
    \caption{Mirroring space}
    \label{fig:mirrored_domain}
\end{figure*}

\section{Algorithm adjustment for quantum computation}

In this section an Hybrid adaptation of the algorithms is presented in order to develop a QFT based Homogenization algorithm. As a starting point, a presentation of the Quantum notations for Quantum computing is given. This will be followed by a description of the Quantum Fourier Transform (QFT) and some useful extensions such as multi-dimensional QFT and QFT on tensor fields. To continue, 2 methods of Quantum Fourier coefficients estimation will be presented. Lastly, a set of results for the different methods will be showed.
\subsection{Quantum notations}

In the following sections, Dirac's bra-ket notation \cite{dirac1939new} will be use to describe quantum states and quantum operations. A quantum superposed state of $n$ qubits $\ket{q}{}_n \in \qH^{\mal n}$ is expressed as followed : $\ket{q}{}_n= \sum_{k=0}^{2^n - 1}x_k\ket{k}{}_n$ where $\ket{k}{}_n$ is called a pure state and is equivalent to the $k-$th vector of the standard basis of $\C^{2^n}$ and where $\sum_{k=0}^{2^n - 1}|x_k|^2=1,\ \boldsymbol{x} =(x_k)_{k \in [0,2^n-1]} \in \C^{2^n}$. A quantum gate is a unitary operation on $\qH^{\mal n}$,then, any quantum operation $\boldsymbol{U}$ is in $\mathbb{U}(\qH^{\mal n})$ and so verify $\boldsymbol{UU}^{\dagger}=\boldsymbol{U}^{\dagger}\boldsymbol{U}=\boldsymbol{I}$ where $\boldsymbol{U}^{\dagger}$ is the self adjoint of $\boldsymbol{U}$ and $\boldsymbol{I}$ is the Identity matrix. For example in the following section three 1 qubit gates will be needed : $$\begin{aligned}&\boldsymbol{H}=\frac{1}{\sqrt{2}} \begin{pmatrix} 1 & 1\\ 1 & -1 \end{pmatrix},  \quad \boldsymbol{S}^{\dagger}= \begin{pmatrix} 1 & 0\\ 0 & -i \end{pmatrix}, \quad \\ & \boldsymbol{U_1}(\lambda)= \begin{pmatrix} 1 & 0\\ 0 & e^{i\lambda} \end{pmatrix}, \quad \forall \lambda \in \R,\end{aligned}$$ where $i$ stands for the imaginary unit defined by $i^2 = 1$. For multi-qubits gates, controlled operations are used. For example, set $\boldsymbol{U} \in \mathbb{U}(\qH)$ a 1 qubit unitary operation (or gate), then the associated Controlled-$\boldsymbol{U}$ operation ($\boldsymbol{C_U}$) acting on 2 qubits is defined by $$\boldsymbol{C_U}= \begin{pmatrix} \boldsymbol{I} & \boldsymbol{0}_2\\ \boldsymbol{0}_2 & \boldsymbol{U} \end{pmatrix}.$$ 
This process can be generalized to $n-$qubits multi-controlled gates. As a quantum algorithm is a sequence of quantum gates, it is a unitary operation on $\qH^{\mal n}$.
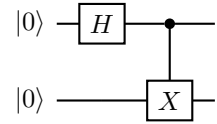
\begin{figure}[H]
    \centering
    \begin{tikzpicture}
\node[scale=1] {

\begin{quantikz}[column sep=0.3cm]
\lstick{$\ket{0}$} & \gate{H} &   \ctrl{+1}  & \qw \\
\lstick{$\ket{0}$} & \qw & \gate{X}  & \qw \\
\end{quantikz}
};
\end{tikzpicture}
    \caption{2 qubits Bell state generation circuit}
    \label{fig:bell}
\end{figure}
As a little example the circuit of 2 qubits Bell state generation is plotted in Figure \ref{fig:bell}. Starting from the initial state $\ket{\psi}=\ket{0}\ket{0}$ an $\boldsymbol{H}$  gate is applied on the first qubit giving the following state : $\boldsymbol{H}\ket{\psi}=\frac{1}{\sqrt{2}}(\ket{1}{}+\ket{0}{})\ket{0}{} $ then a $\boldsymbol{C_X}$ gate is applied using the first qubit as controller : $\boldsymbol{C_XH}\ket{\psi}=\frac{1}{\sqrt{2}}(\ket{11}{}+\ket{00}{})$. The state produced is called an EPR state (or a Bell state) \cite{nielsen2002quantum}. It can be noticed that measuring one qubit will determine the outcome of the measurement of the second one, that phenomenon is called entanglement.

\subsection{Amplitude encoding}

The first main challenge of Quantum algorithm is usually in the initial state encoding. This can be summarized in an operator $\boldsymbol{\mathcal{I}}_{\boldsymbol{f}} \in \mathbb{U}(\qH^{\mal n})$, it maps a normalized vector $\boldsymbol{f}$ of $N=2^n$ coefficients with sum equal to $1$ following \begin{equation}
\label{eqn:mapping}
    \boldsymbol{\mathcal{I}}_{\boldsymbol{f}}\ket{0}{}_n=\sum_{k=0}^{N - 1}f_k\ket{k}{}_n.
\end{equation}
This operation can be implemented using uniformly controlled rotations as explained by Mottonen and al. in \cite{mottonen2004transformation}. Every controlled rotations can be represented with the following operators which are directly available in Qiskit, $$\begin{aligned}&\boldsymbol{R_x}(\theta)= \begin{pmatrix} \cos(\theta/2) & -i\sin(\theta/2)\\ -i\sin(\theta/2) & \cos(\theta/2) \end{pmatrix}, \\&  \boldsymbol{R_y}= \begin{pmatrix} \cos(\theta/2) & -\sin(\theta/2)\\ \sin(\theta/2) & \cos(\theta/2) \end{pmatrix}, \\ & \boldsymbol{R_z}(\theta)= \begin{pmatrix} e^{-i\theta/2} & 0\\ 0 & e^{i\theta/2} \end{pmatrix},\quad \forall \theta \in \R.\end{aligned}$$
According to Mottonen and al., the operator $\boldsymbol{\mathcal{I}}_{\boldsymbol{f}}$ needs in total $2^{n+2}-4n-4$ $\boldsymbol{C_X}$ gates and $2^{n+5}-5$ one-qubit elementary gates. Therefore, the complexity of the encoding is $\mathcal{O}(N)$.

\subsection{Quantum Fourier Transform}

\begin{figure*}
\centering
\begin{tikzpicture}
\node[scale=0.7] {
\begin{quantikz}[column sep=0.3cm]
\lstick{$\ket{q_1}$} & \qw & \gate[5]{\text{QFT}_n} & \qw \\
\lstick{$\ket{q_2}$} & \qw & & \qw \\ 
\vdots &&&& \\
\lstick{$\ket{q_{n-1}}$}& \qw & & \qw  \\
\lstick{$\ket{q_{n}}$}& \qw & & \qw
\end{quantikz}
=
\begin{quantikz}[column sep=0.3cm]
\lstick{$\ket{q_1}$} & \gate{H} &   \gate{U_1(\frac{\pi}{2})} & \ \dots \ \qw & \gate{U_1(\frac{\pi}{2^{n-2}})} & \gate{U_1(\frac{\pi}{2^{n-1}})} & \qw & \qw & \qw & \qw & \qw & \qw & \qw & \qw &  \swap{4} & \qw & \qw \\
\lstick{$\ket{q_2}$} & \qw & \ctrl{-1} & \qw & \qw & \qw & \gate{H} & \ \dots \ \qw &  \gate{U_1(\frac{\pi}{2^{n-3}})} & \gate{U_1(\frac{\pi}{2^{n-2}})} & \qw & \qw & \qw & \qw & \qw & \swap{2} & \qw \\
\vdots &&&&&&&&&&&&&&&& \ \ldots \ \\ 
\lstick{$\ket{q_{n-1}}$}& \qw & \qw &  \qw  & \ctrl{-3} & \qw & \qw & \qw & \ctrl{-2} & \qw & \gate{H} & \qw &  \gate{U_1(\frac{\pi}{2})} & \qw & \qw & \targX{} &  \qw\\
\lstick{$\ket{q_{n}}$} &\qw & \qw &  \qw  & \qw& \ctrl{-4} & \qw & \qw & \qw  & \ctrl{-3} & \qw & \qw & \ctrl{-1} & \gate{H} & \targX{} & \qw & \qw
\end{quantikz}
};
\end{tikzpicture}
    \caption{QFT Schematic }
    \label{fig:QFT}
\end{figure*}
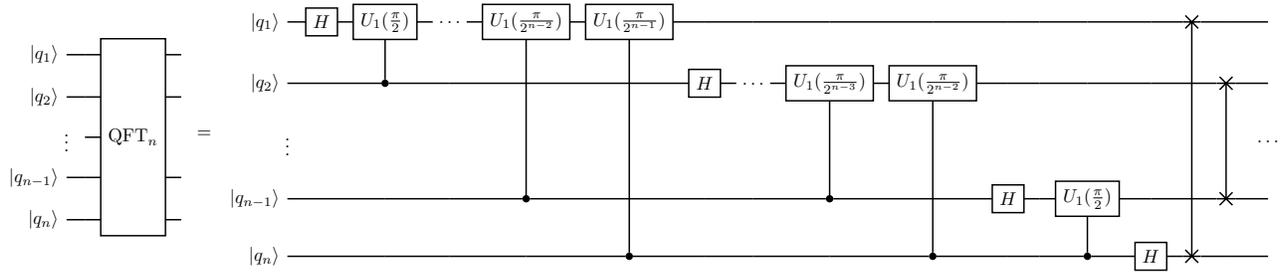

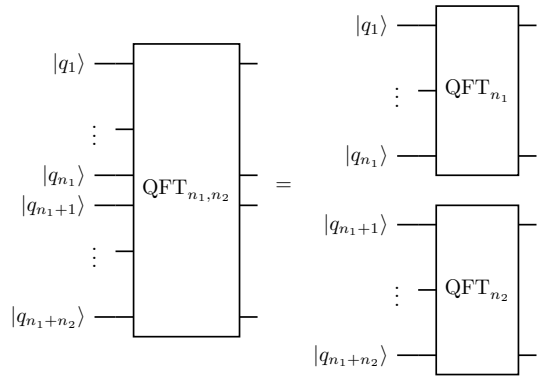
\begin{figure}[H]
\centering
\begin{tikzpicture}
\node[scale=0.8] {

\begin{quantikz}[column sep=0.3cm]
\lstick{$\ket{q_1}$} & \qw & \gate[6]{\text{QFT}_{n_1,n_2}} & \qw \\
\vdots &&&& \\
\lstick{$\ket{q_{n_1}}$} & \qw & & \qw \\ 
\lstick{$\ket{q_{n_1+1}}$}& \qw & & \qw  \\
\vdots &&&& \\
\lstick{$\ket{q_{n_1+n_2}}$}& \qw & & \qw
\end{quantikz}
=
\begin{quantikz}[column sep=0.3cm]
\lstick{$\ket{q_1}$} & \qw & \gate[3]{\text{QFT}_{n_1}} & \qw \\
\vdots &&&& \\
\lstick{$\ket{q_{n_1}}$} & \qw & & \qw \\ 
\lstick{$\ket{q_{n_1+1}}$}& \qw & \gate[3]{\text{QFT}_{n_2}}& \qw  \\
\vdots &&&& \\
\lstick{$\ket{q_{n_1+n_2}}$}& \qw & & \qw
\end{quantikz}
};
\end{tikzpicture}
    \caption{2 dimensions QFT Schematic}
    \label{fig:QFT_2D}
\end{figure}

Used as a building block of numerous Quantum algorithms, the Quantum Fourier Transform (QFT) \cite{jozsa1998quantum} is a well documented algorithm.  The aim is to perform the Normalized Discrete Fourier transform of normalized input encoded in the different coefficients of a quantum superposed state. For a given input normalized vector $\boldsymbol{f}$ of size $N=2^n$, the number of qubit needed to encode the vector is $n$ and the Quantum Fourier transformed is expressed as followed, \begin{equation}
    \resizebox{0.9 \hsize}{!}{$\begin{split}
         \text{QFT}_n : \  & \hspace{3.6em} \qH^{\mal n}  \hspace{1.9em} \longrightarrow  \hspace{3.9em} \qH^{\mal n} \\
        & \ket{q}{}_n= \sum_{k=0}^{N - 1}f_k\ket{k}{} \longmapsto \ \ket{\hat{q}}{}_n= \sum_{k=0}^{N - 1}\hat{f}_k\ket{k}{},
    \end{split}$}
\end{equation} where $\hat{f}_k=\frac{1}{\sqrt{N}} \sum_{j=0}^{N - 1} \omega_N^{jk}f_j$, and $\omega_N=e^{-\frac{2i \pi}{N}}$. To implement this unitary operation we use the gate sequence described on Figure \ref{fig:QFT}. The complexity of the algorithm can be estimated by the number of quantum gates used which here is in $\mathcal{O}(n^2)=\mathcal{O}(\text{log}(N)^2)$.
This operation can easily be extented in multi-dimension. For example in 2 Dimensions with $N=N_1 \times N_2 = 2^{n_1} \times 2^{n_2}$ , two quantum registers of respectively $n_1$ and $n_2$ qubits are used, \begin{equation}
         \ket{\hat{q}}_{n_1+n_2}= \sum_{j=0}^{N_1 - 1}\sum_{k=0}^{N_2 - 1}\hat{f}_{jk}\ket{j}{}_{n_1}\ket{k}{}_{n_2},
\end{equation}  where $\hat{f}_{jk}= \frac{1}{\sqrt{N_1 N_2}}\sum_{l=0}^{N_1 - 1}\sum_{m=0}^{N_2 - 1}\omega_{N_1}^{jl}\omega_{N_2}^{km} f_{lm}$. This is equivalent of applying a 1D QFT on both registers simultaneously as you can see on Figure \ref{fig:QFT_2D}. Similarly, the complexity of this algorithm is determined by the number of quantum gates which here is in $\mathcal{O}(n_1^2+n_2^2)=\mathcal{O}(\text{log}(N_1)^2+\text{log}(N_2)^2)$. Nevertheless the two operations are independant, meaning that they can be parallelized reducing the complexity to $\mathcal{O}(\text{max}(n_1^2,n_2^2))=\mathcal{O}(\text{max}(\log(N_1)^2,\log(N_2)^2))$

As a first optimization of the QFT, a method from \cite{10.5555/1972505} can be used. using the fact that the swap gates are only used to get the ordering of the vector right, it can be deleted without any impact on the outcome, the reordering can then happen classically.

\subsection{QFT extension to tensor fields}

\begin{figure}[H]
\centering
\begin{tikzpicture}
\node[scale=0.8] {
\begin{quantikz}[column sep=0.3cm]
\lstick{$\ket{q_1}$} & \qw & \gate[6]{\text{QFT}_{n_1,n_2}} & \qw \\
\vdots &&&& \\
\lstick{$\ket{q_{n_1}}$} & \qw & & \qw \\ 
\lstick{$\ket{q_{n_1+1}}$}& \qw & & \qw  \\
\vdots &&&& \\
\lstick{$\ket{q_{n_1+n_2}}$}& \qw & & \qw \\
\lstick{$\ket{n_1+n_2+1}$} &\qw&\qw& \qw& \\
\vdots &&&& \\
\lstick{$\ket{n_1+n_2+n_d}$} &\qw&\qw& \qw& \\
\end{quantikz}
};
\end{tikzpicture}
    \caption{2 dimensions QFT Schematic on tensor field}
    \label{fig:QFT_tensor}
\end{figure}
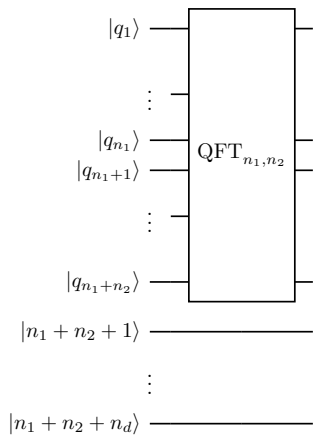

Another good advantage of Quantum Fourier Transform is that it can easily be extended on tensor fields vectors. For example, suppose that we have a discrete tensor of order 2 $\boldsymbol{\nu}\in \R^{N_1 \times N_2 \times d}$ where $N_1=2^{n_1}$ and $N_2=2^{n_2}$ stand for the number of discrete points in the two dimensional space and $d$ is the number of component of the tensor, it should be of the form $d=2^{n_d}$. Then the tensor field can be implemented in a quantum state $$\ket{\boldsymbol{\nu}}{}_{n_1,n_2,n_d}=\frac{1}{\|\boldsymbol{\nu}\|}\begin{pmatrix} \boldsymbol{\nu}_1 \\ \vdots\\ \boldsymbol{\nu}_d
\end{pmatrix}.$$ Which is implemented by adding another quantum register of $n_d$ qubits to the spatial register of $n_1+n_2$ qubits : 
\begin{equation}
        \ket{\boldsymbol{\nu}}_{n_1,n_2,n_d}= \sum_{k=1}^{d}\ket{\boldsymbol{\nu}_k}_{n_1,n_2}\ket{k}{}_{n_d}
\end{equation} where $\ket{\boldsymbol{\nu}_k}_{n_1,n_2}=\frac{\boldsymbol{\nu_k}}{\|\boldsymbol{\nu}\|}$ is the $k-\text{th}$ component of $\ket{\boldsymbol{\nu}}{}_{n_1,n_2,n_d}$. It can be noticed that \begin{equation}
        \ket{\hat{\boldsymbol{\nu}}}_{n_1,n_2,n_d}= \sum_{k=1}^{d} \text{QFT}_{n_1,n_2}\ket{\boldsymbol{\nu}_k}_{n_1,n_2}\ket{k}{}_{n_d},
\end{equation}the 2D-QFT can be performed by applying the 2D-QFT on the dimensional register of $n_1+n_2$ qubit (Figure \ref{fig:QFT_tensor}).
In case the number of components of the tensor field is not a power of $2$ the algorithm still works, the non-existing components just have to be set to $\boldsymbol{0}$, for example suppose that $\boldsymbol{\nu}\in \R^{N_1 \times N_2 \times d-1}$ with $d=2^n_d$, then $$\ket{\boldsymbol{\nu}}{}_{n_1,n_2,n_d}=\frac{1}{\|\boldsymbol{\nu}\|}\begin{pmatrix} \boldsymbol{\nu}_1 \\ \vdots\\ \boldsymbol{\nu}_{d-1} \\ \boldsymbol{0}
\end{pmatrix}.$$

\subsection{Quantum coefficient estimation using Hadamard Test }

To replace the classical FFT by a QFT, one of the challenge is to determine the different coefficients of the superposed quantum state. More specifically, the aim is for a given input state $\ket{\psi}{}_n \in \qH^{\otimes n}$ to estimate $\braket{\psi}{k}$ and this $\forall k\in[0,2^n -1]$. To do so a method can be deduced from the Hadamard test. Hadamard test method is made to estimate for a given input state $\ket{\psi}{}_n \in \qH^{\otimes n}$ and a unitary operator $\boldsymbol{U}\in \mathbb{U}(\qH^{\otimes n})$ : $\bra{\psi}{}\boldsymbol{U}\ket{\psi}{}$ by adding an ancilla qubit to the quantum register where $\ket{\psi}{}$ is encoded. In fact, using the circuit described on Figure \ref{fig:ReHadamard} leads to the following property, \begin{equation}
    2P(\ket{1}{})-1=\text{Re}(\bra{\psi}{}\boldsymbol{U}\ket{\psi}{}),
\end{equation} where $P(\ket{1}{})$ stands for the probability of the ancilla qubit to output one after measurement. This makes it possible to estimate the real part of $\bra{\psi}{}\boldsymbol{U}\ket{\psi}{}$ : $\text{Re}(\bra{\psi}{}\boldsymbol{U}\ket{\psi}{})$.  Respectively, using the circuit described on Figure \ref{fig:IMHadamard} gives a similar property for the imaginary part, \begin{equation}
    2P(\ket{1}{})-1=\text{Im}(\bra{\psi}{}\boldsymbol{U}\ket{\psi}{}).
\end{equation} In order to measure each state coefficients what is needed is $\braket{\psi}{k} \forall k \in [0,2^n-1]$. The operator U can then be transformed in order to have $\braket{\psi}{k} $. Let define $\boldsymbol{U}_{\psi\xrightarrow{}0} \in \mathbb{U}(\qH^{\otimes n}) $ such as for a given state $\ket{\psi}, \ \boldsymbol{U}_{\psi\xrightarrow{}0}\ket{\psi}{}_n=\ket{0}{}_n$ and $\boldsymbol{U}_{0\xrightarrow{}k} \in \mathbb{U}(\qH^{\otimes n}) $ such as $\boldsymbol{U}_{0\xrightarrow{}k}\ket{0}{}_n=\ket{k}{}_n$. Then the operator $\boldsymbol{U}=\boldsymbol{U}_{0\xrightarrow{}k} \boldsymbol{U}_{\psi\xrightarrow{}0}$ leads for the circuit shown in Figure \ref{fig:ReHadamard} to \begin{equation}
    2P(\ket{1}{})-1=\text{Re}(\braket{\psi}{k}),
\end{equation} and for the circuit shown in Figure  \ref{fig:IMHadamard} to : \begin{equation}
    2P(\ket{1}{})-1=\text{Im}(\braket{\psi}{k}).
\end{equation} The operator $\boldsymbol{U}_{0\xrightarrow{}k}$ is trivial. The operator $\boldsymbol{U}_{\psi\xrightarrow{}0}$ is a bit more complicated. The circuit to encode the state $\ket{\psi}{}_n$ is in general known and defined by $\boldsymbol{\mathcal{I}}_{\psi}$. The operator $\boldsymbol{U}_{\psi\xrightarrow{}0}$ can be defined as $\boldsymbol{U}_{\psi\xrightarrow{}0}=\boldsymbol{\mathcal{I}}_{\psi}^{\dagger}$. This operation have a circuit easy to build as all the gate composing the initialization circuit have trivial inverse gate.
\begin{figure}[H]
     \centering
     \begin{subfigure}[b]{0.4\textwidth}
\centering
\begin{tikzpicture}
\node[scale=1] {

\begin{quantikz}[column sep=0.3cm]
\lstick{$\ket{0}$} & \gate{H} &   \ctrl{+1} &  \gate{H}  & \meter{} &  \\
\lstick{$\ket{\psi}_n$} & \qw & \gate{U} & \qw & \qw & \\
\end{quantikz}
};
\end{tikzpicture}
    \caption{Real part Hadamard test circuit}
    \label{fig:ReHadamard}
     \end{subfigure}
     \hfill
     \begin{subfigure}[b]{0.4\textwidth}
\centering
\begin{tikzpicture}
\node[scale=1] {

\begin{quantikz}[column sep=0.3cm]
\lstick{$\ket{0}$} & \gate{H} & \gate{S^{\dagger}} & \ctrl{+1} &  \gate{H}  & \meter{} &  \\
\lstick{$\ket{\psi}_n$} & \qw & \qw & \gate{U} & \qw & \qw & \\
\end{quantikz}
};
\end{tikzpicture}
    \caption{Imaginary part Hadamard test circuit}
    \label{fig:IMHadamard}
     \end{subfigure}
        \caption{Hadamard test circuits}
        \label{fig:Hadamard_circuits}
\end{figure}
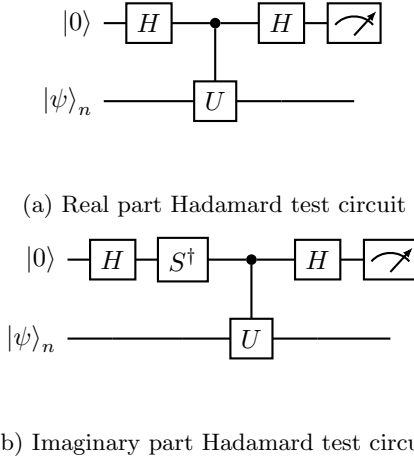 The method can now be adapted to the determination of the Quantum Fourier transform coefficients. From the initialization operator $\boldsymbol{I}_f$ the operator $\boldsymbol{U}_{\psi\xrightarrow{}0}$ can be defined : \begin{equation}
    \boldsymbol{U}_{\psi\xrightarrow{}0}=\boldsymbol{\mathcal{I}}_{\boldsymbol{f}}^{\dagger}\mathcal{QFT}_n^{\dagger}.
\end{equation}
In the next section different methods to estimate $P(\ket{1}{})$ with a given accuracy will be described. Finally, let's define $\mathcal{H}_{\R}(\boldsymbol{U}) : \mathbb{U}(\qH^{\mal n}) \xrightarrow{} \mathbb{U}(\qH^{\mal n})$ which from a given operator $\boldsymbol{U}$ returns the associate real Hadamard test circuit. Respectively, let's define $\mathcal{H}_{i\R}(\boldsymbol{U}) : \mathbb{U}(\qH^{\mal n}) \xrightarrow{} \mathbb{U}(\qH^{\mal n})$ which from a given operator $\boldsymbol{U}$ returns the associate Imaginary Hadamard test circuit.

\subsection{Amplitude amplification}
 As the exact determination of the Quantum coefficients is not possible in practice, different methods have been developed in order to estimate the amplitude of a quantum coefficient. Most of these methods are based on quantum amplification developed by Brassard \& al. \cite{brassard2002quantum} as a generalization of the Grover algorithm \cite{grover1996fast}. Suppose that the amplitude $x_k$ of base state $\ket{k}{}_n$ has to be estimated. Set the operator $\boldsymbol{\mathcal{A}}_k \in \mathbb{U}(\qH^{\otimes n})$ where $k\in[0,N-1], \ N=2^n$ acting for a given vector $\boldsymbol{x} =(x_k)_{k \in [0,2^n-1]} \in \qH^{\mal n}$ as, \begin{equation}
\label{eq:A_i}
    \ket{q}{}_{n+1}=\boldsymbol{\mathcal{A}}_k(\boldsymbol{x})\ket{0}{}_n\ket{0}{}= x_k\ket{k}_n{}\ket{1}{} + (\sum_{l \neq k} x_l\ket{l}{}_n)\ket{0}{}.
\end{equation} This will first, initialize the state $\ket{x}{}_n = \sum_{k=0}^{2^n-1}x_k\ket{k}{}_n$. Then it will be followed by the entanglement of the base state $\ket{k}{}_n$ to the state $\ket{1}{}$ of an ancilla qubit. As $\ket{q}{}_{n+1} \in \qH^{\mal n}$ there exist $\theta_k \in [0, \frac{\pi}{2}]$ satisfying $|x_k|=\cos{\theta_k}$ and $\sqrt{\sum_{l \neq k}|x_l|^2}= \sin{\theta_k}$. 
Set now Grover's operator introduced in \cite{grover1996fast}, \begin{equation}
    \boldsymbol{\mathcal{Q}}_k(x)= \boldsymbol{\mathcal{A}}_k(x) \boldsymbol{\mathcal{S}}_0 \boldsymbol{\mathcal{A}}_k^{-1}(x) \boldsymbol{\mathcal{S}}_{k},
\end{equation} where $\boldsymbol{\mathcal{S}}_{k}$ changes the amplitude phase of the state $k$, 
\begin{equation}
\label{eqn:chi}
    \ket{x}{} \rightarrow{}
    \begin{cases} -\ket{x}{} \text{ if } x= k, \\ \ket{x}{} \text{ otherwise }. \end{cases}
\end{equation} In the end the application of Grover's operator will lead to an amplification such as if the operator is applied $m$ times, \begin{equation}\begin{aligned}
        \boldsymbol{\mathcal{Q}}_k(\boldsymbol{x})^m\ket{q}{}_{n+1} = & \cos\big((2m+1)\theta_k\big)e^{\phi_k}\ket{k}_n{}\ket{1}{}+ \\ & \sin\big((2m+1)\theta_k\big)(\sum_{l \neq k}e^{\phi_l}\ket{l}{}_n)\ket{0}{},
        \end{aligned}
    \end{equation}  where $\phi_l, l\in [0,N-1]$ are the phases of the different coefficients $x_l$. This amplification is a well known result that is explained in \cite{brassard2002quantum}.

\subsection{Amplitude estimation}

Most of current amplitude estimation methods are based on the amplitude amplification. The usual way of doing it is to determine the eigen value of Grover's operator. To do so, it can be noticed that the amplitude $x_k=cos((2m+1)\theta_k)^2$ is a periodic function in $m$ with period $2\theta_k$. The aim of current amplitude estimation methods is to determine this $\theta_k$. For example, the method developed by Brassard \& al. in \cite{brassard2002quantum} is using a register of multiple ancilla qubits, and by using the phase estimation algorithm of Nielsen and al. \cite{nielsen2002quantum} it estimates $\theta_k$. It can be noticed that Brassard's method has an asymptotic convergence following $|\Tilde{x}_k-x_k|=\mathcal{O}(\frac{1}{2^m})$ where $m$ is the size of the ancilla register. This limitation of a need of a big ancilla register to have an accurate estimation makes the method hard to use in practice due to hardware apparent noise. By contrast, recent methods are trying to limit them-self to only one ancilla qubit that will be measured multiple times with different number of amplifications, avoiding the use of the phase estimation method. This is the case of a method developed by Suzuki and al. \cite{suzuki2020amplitude} called Most Likelihood Quantum Amplitude Estimation (MLQAE) that will be used in the following sections. This method is using most likelihood functions based on the results of quantum amplification experiments with various number of amplifications to estimate the periodicity of Grover's operator. A set of amplification experiments can be define with a set of amplification numbers $\boldsymbol{M}=\{m_1, \dots ,m_M \}$ and a number of repetitions $S$ (more comonly called number of shots). An experiment consists in the application of $ \boldsymbol{\mathcal{Q}}_k(\boldsymbol{x})^m{\boldsymbol{\mathcal{A}}_k(\boldsymbol{x})\ket{0}{}_n\ket{0}{}}$ wit $m$ taken in $\boldsymbol{M}$, then, the ancilla qubit is measured. The experiments are repeated $S$ time for each values of amplifications number in $\boldsymbol{M}$. In the end the results of this experiments are used to create a function $\boldsymbol{L}(\theta)$ with a maximum value obtained for the most likelihood value of $\theta_k$ fitting with the experiment results. The original paper of Suzuki and al. \cite{suzuki2020amplitude} shows that the optimal set of amplification numbers is $\boldsymbol{M}=\{0\}\cup\{1,2, \dots 2^{M-1}\}$ leading to the following asymptotic comportment of the error, \begin{equation}\label{eq:error_mlqae}e(\boldsymbol{M})=|\Tilde{x}_k-x_k|=\mathcal{O}(\frac{1}{2^{M-1}\sqrt{S}}).\end{equation} The method complexity depends on the complexity of the operator $\boldsymbol{\mathcal{A}}_k(\boldsymbol{x})$ and on the targeted error $e$. Setting $\boldsymbol{\mathcal{A}}_k(\boldsymbol{x})$ complexity to C($\boldsymbol{\mathcal{A}}_k(\boldsymbol{x})$) leads to a complexity on quantum computer in $\mathcal{O}(e^{-1}\boldsymbol{\mathcal{A}}_k(\boldsymbol{x}))$. It is also important to notice that the algorithm has a classical component which has a complexity only depending on the error and acting in $\mathcal{O}(e^{-1}\text{ln} (e^{-1}))$.

\subsection{Application to QFT coefficients estimation}

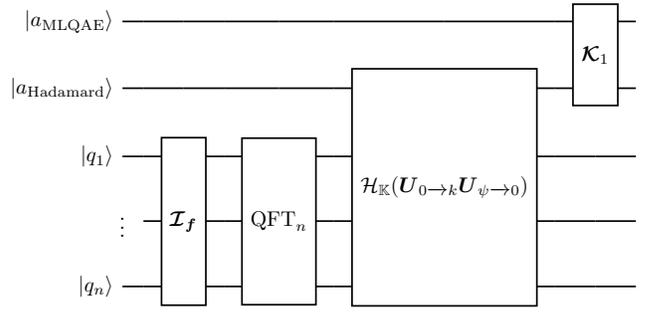
\begin{figure}
\centering
\begin{tikzpicture}
\node[scale=0.8] {
\begin{quantikz}[column sep=0.3cm]
\lstick{$\ket{a_{\text{MLQAE}}}$} & \qw &  \qw &  \qw &  \qw &  \qw &\qw & \qw & \gate[2]{\boldsymbol{\mathcal{K}}_1} & \qw \\
\lstick{$\ket{a_{\text{Hadamard}}}$} & \qw &  \qw & \qw &  \qw & \qw  &  \gate[4]{\mathcal{H}_{\mathbb{K}}(\boldsymbol{U}_{0\xrightarrow{}k}\boldsymbol{U}_{\psi\xrightarrow{}0})} & \qw & & \qw \\ 
\lstick{$\ket{q_1}$}& \qw & \gate[3]{\boldsymbol{\mathcal{I}}_{\boldsymbol{f}}} & \qw & \gate[3]{\text{QFT}_n} & \qw & & \qw & \qw & \qw \\
\vdots &  & & \qw & & \qw & & \qw & \qw & \qw  \\
\lstick{$\ket{q_{n}}$}& \qw & & \qw & & \qw & & \qw & \qw & \qw
\end{quantikz}
};
\end{tikzpicture}
    \caption{$\boldsymbol{\mathcal{A}}_k$ operator for Hadamard test}
    \label{fig:A_i}
\end{figure}

As MLQAE is based on Grover's amplitude amplification algorithm, the only thing needed to arrive to the estimation of the Fourier coefficient stored in the coefficient of the base state $\ket{k}_n$ is to transform the operator $\boldsymbol{\mathcal{A}}_k$. To do so, let's first define $\boldsymbol{\mathcal{K}}_k \in \mathbb{U}(\qH^{\mal n+1})$ an operator acting on an arbitrary state $\ket{q}{}_n=\sum_{k=0}^{N-1}x_k\ket{k}{}_n$ which will entangled the base state $\ket{k}_n$ to the state \ket{1}{} of an ancilla qubit, \begin{equation}
    \boldsymbol{\mathcal{K}}_k\ket{q}{}_n ,\ket{0}{}=x_k\ket{k}{}_n\ket{1}{}+(\sum_{l\neq i}x_l\ket{l}{})\ket{0}{}. 
\end{equation} It can be seen as the operator $\boldsymbol{\mathcal{A}}_k$ without the initialization part. In the case of the Fourier transform this "initialization" operation could be decomposed into the initialization of the input state $\boldsymbol{f}$ followed by the QFT, \begin{equation}
    \boldsymbol{\mathcal{A}}_k=\boldsymbol{\mathcal{K}}_k\text{QFT}_n\boldsymbol{\mathcal{I}}_{\boldsymbol{f}}.
\end{equation} Nevertheless this would lead to the estimation of the amplitude of the coefficient loosing the information about the complex phase. To solve this problem we have to add the application of the Hadamard circuits to the operator $\boldsymbol{\mathcal{A}}_k$ to estimate first the real part of the coefficient, and then the imaginary part. This leads to the addition of another ancilla qubit for the Hadamard test. Now the amplitude of the $\ket{1}{}$ state of the Hadamard ancilla qubit has to be estimated. To do so, $\boldsymbol{\mathcal{K}}_1$ is used with $\ket{a_{\text{Hadamard}}}{}$ as target register and $\ket{a_{\text{MLQAE}}}{}$ as ancilla qubit. Recalling that in the case of the QFT $\boldsymbol{U}_{\psi\xrightarrow{}0}=\boldsymbol{\mathcal{I}}_{\boldsymbol{f}}^{\dagger}\text{QFT}_n^{\dagger}$ the new definition of $\boldsymbol{\mathcal{A}}_k$ is then, \begin{equation}
    \boldsymbol{\mathcal{A}}_k=\boldsymbol{\mathcal{K}}_1\mathcal{H}_{\mathbb{K}}(\boldsymbol{U}_{0\xrightarrow{}k}\boldsymbol{U}_{\psi\xrightarrow{}0})\text{QFT}_n\boldsymbol{\mathcal{I}}_{\boldsymbol{f}},
\end{equation} where $\mathcal{H}_{\mathbb{K}}$ stands for the operator of Hadamard test estimating the real part ($\mathcal{H}_{\mathbb{R}}$) or on the imaginary part of the numbers ($\mathcal{H}_{i\mathbb{R}}$). The circuit of this operator is shown on Figure \ref{fig:A_i}. The complexity of this operation $C(\boldsymbol{\mathcal{A}}_k)$ is in $O(N)$ as the limiting operation remains the initialization leading to a global complexity on quantum computers in $\mathcal{O}(e^{-1}N)$.

To determine the Discrete Fourier transform ($\text{DFT}$) of an arbitrary vector $\boldsymbol{f}=(f_k)_{k \in [0,N-1]}$ using Quantum Fourier transform and MLQAE estimation, a re-normalization has to be considered as $\text{QFT}_n$ can only be used on unitary vector : \begin{equation}
    \text{DFT}(\boldsymbol{f})=
    \begin{cases}
    \boldsymbol{0} \text{  if } \|\boldsymbol{f}\|=0, \\
    \sqrt{N}\|\boldsymbol{f}\|\text{QFT}_n(\frac{\boldsymbol{f}}{\|\boldsymbol{f}\|}) \text{   otherwise}.
    \end{cases}
\end{equation}
An algorithm of QFT based Fourier coefficients estimation can then be deduced.

\begin{algorithm}[H]
\caption{Hadamard QFT}\label{alg:Hadamard_QFT_based_DFT}
\begin{algorithmic}
\Require $\boldsymbol{f},\boldsymbol{M}$
\For{$k=1,N$}
\State $\boldsymbol{\mathcal{A}}_k=\boldsymbol{\mathcal{K}}_1\mathcal{H}_{\mathbb{R}}(\boldsymbol{U}_{0\xrightarrow{}k}\boldsymbol{U}_{\psi\xrightarrow{}0})\text{QFT}_n\boldsymbol{\mathcal{I}}_{\boldsymbol{f}},$
\State $\boldsymbol{\mathcal{Q}}_k= \boldsymbol{\mathcal{A}}_k \boldsymbol{\mathcal{S}}_0 \boldsymbol{\mathcal{A}}_k^{-1}\boldsymbol{\mathcal{S}}_{\chi},$
\State $a_k=\sqrt{N} \|f\| \mlqae (\boldsymbol{\mathcal{Q}}_k,\boldsymbol{M}),$
\EndFor
\For{$k=1,N$}
\State $\boldsymbol{\mathcal{A}}_k=\boldsymbol{\mathcal{K}}_1\mathcal{H}_{i\mathbb{R}}(\boldsymbol{U}_{0\xrightarrow{}k}\boldsymbol{U}_{\psi\xrightarrow{}0})\text{QFT}_n\boldsymbol{\mathcal{I}}_{\boldsymbol{f}},$
\State $\boldsymbol{\mathcal{Q}}_k= \boldsymbol{\mathcal{A}}_k \boldsymbol{\mathcal{S}}_0 \boldsymbol{\mathcal{A}}_k^{-1}\boldsymbol{\mathcal{S}}_{\chi},$
\State $b_k=\sqrt{N} \|\boldsymbol{f}\| \mlqae (\boldsymbol{\mathcal{Q}}_k,\boldsymbol{M}),$
\EndFor
\return $\boldsymbol{a}+i\boldsymbol{b}.$
\end{algorithmic}
\end{algorithm} The error of approximation for a given coefficient can now be defined as the error of MLQAE method (Equation \eqref{eq:error_mlqae}) multiplied by the norm of the Fourier coefficients vector, \begin{equation}
    e_k(M)=|\Tilde{a}_k+i\Tilde{b}_k-a_k-ib_k|=O(\frac{\sqrt{N}\|\boldsymbol{f}\|}{2^{M-1}\sqrt{S}}).
    \label{eqn:asymptotic_error}
\end{equation}

\subsection{Efficient estimation method for symmetric input domain}

In this section, a method to improve the Hadamard test method for specific boundary conditions is developed in order to avoid big controlled operations. In fact, in the case of PMBC on mirrored space it can be noticed that as the input data domain is symmetric or anti-symmetric, its Fourier transform verify certain properties,  following Grimm-Strele and Kabel method \cite{GrimmStrele2021FFTB}, for a vector $\boldsymbol{f}=(f_k)_{k  \in [0,N-1]} \in \mathbb{R}^N$ where $f_j = f_{N-j-1}$ for even symmetry and $f_j = -f_{N-j-1}$ for odd symmetry, then : \begin{equation}
    \hat{f_k}=\frac{2e^{i\beta_k}}{\sqrt{N}}\sum_{j=0}^{N/2-1} f_j\text{trig}(\frac{\pi(2j+1)k}{N})
\end{equation} With $\beta_k=\frac{k\pi}{N}$ and $\text{trig}=\cos$ in case of even symmetry and  $\beta_k=-\frac{k\pi}{N}$ and $\text{trig}=\sin$ in case of odd symmetry. This results leads to : \begin{equation}
    \frac{\hat{f_k}}{e^{i\beta_k}} \in \R.
\end{equation} 
Now that it's known that every coefficients $\hat{f}_k$ is in $e^{i\beta_k}\R$, it remains to distinguish whether the coefficient is on the right side or on the left side of the origin (see Figure \ref{fig:graph_phase}).
From there a method in two steps can be deduced using two MLQAE algorithms in a row. The first one will be used to determine the absolute value of the $\hat{f}_j$ coefficients by changing the operator $\boldsymbol{\mathcal{A}}_k$ to : \begin{equation}
    \boldsymbol{\mathcal{A}}_k=\boldsymbol{\mathcal{I}}_{\boldsymbol{f}}\text{QFT}_n\boldsymbol{\mathcal{K}}_k.
\end{equation}

\begin{algorithm}[H]
\caption{QFT amplitude}\label{alg:QFT_based_DFT}
\begin{algorithmic}
\Require $\boldsymbol{f},\boldsymbol{M}$
\For{$k=1,N$}
\State $\boldsymbol{\mathcal{A}}_k=\boldsymbol{\mathcal{I}}_{\frac{\boldsymbol{f}}{\|\boldsymbol{f}\|}}\text{QFT}_n\boldsymbol{\mathcal{K}}_k$
\State $\boldsymbol{\mathcal{Q}}_k= \boldsymbol{\mathcal{A}}_k \boldsymbol{\mathcal{S}}_0 \boldsymbol{\mathcal{A}}_k^{-1}\boldsymbol{\mathcal{S}}_{\chi}$
\State $a_k=\sqrt{N} \|\boldsymbol{f}\| \mlqae (\boldsymbol{\mathcal{Q}}_k,\boldsymbol{M})$
\EndFor
\return $\boldsymbol{a}$
\end{algorithmic}
\end{algorithm}
This part will give an approximation of the amplitudes of the Fourier coefficients.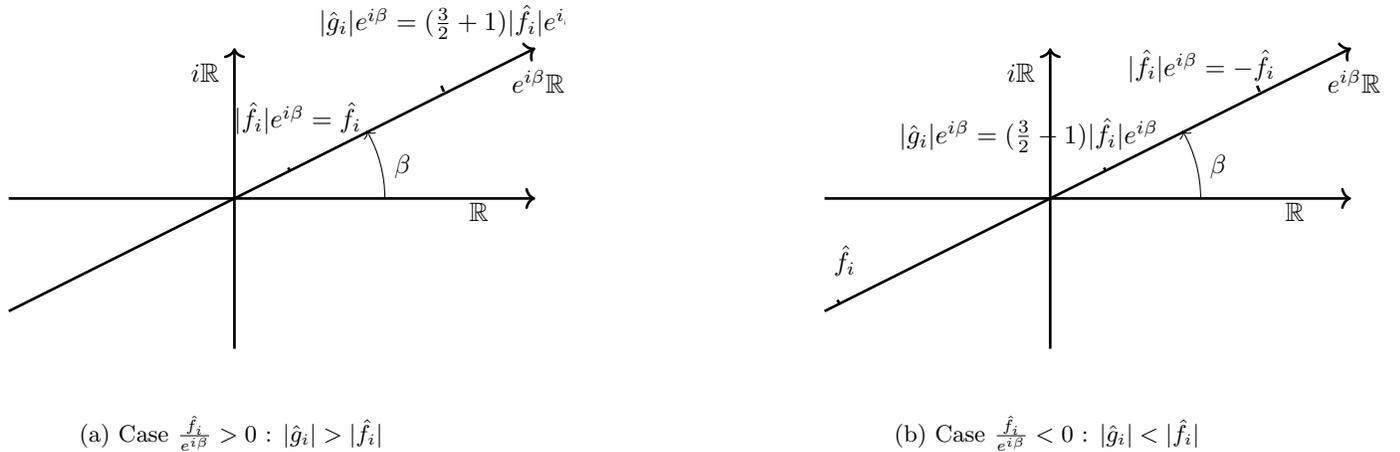
\begin{figure*}
     \centering
     \begin{subfigure}[b]{0.4\textwidth}
     \centering
\begin{tikzpicture}
\clip(-3.6593270712778567,-2.6796850995853734) rectangle (4.392471953927827,2.772429935419367);
\draw [->,line width=1pt] (0,0) -- (4,2);
\draw [line width=1pt] (0,0)-- (-2,-1);
\draw [line width=1pt] (-2,-1)-- (-3,-1.5);
\draw [->,line width=1pt] (0,0) -- (4,0);
\draw [->,line width=1pt] (0,0) -- (0,2);
\draw (3.562017345946972,1.8336551611801402) node[anchor=north west] {$e^{i\beta}\mathbb{R}$};
\draw (3.,0.0734524594815899) node[anchor=north west] {$\mathbb{R}$};
\draw (-0.13,1.4) node[anchor=north west] {$|\hat{f}_i|e^{i\beta}=\hat{f}_i$};
\draw [line width=1pt] (2.75,1.5)-- (2.78474075768685,1.4210640423044958);
\draw (1.,2.7) node[anchor=north west] {$|\hat{g}_i|e^{i\beta}=(\frac{3}{2}+1)|\hat{f}_i|e^{i\beta}$};
\draw (-0.7,1.9239219663954505) node[anchor=north west] {$i\mathbb{R}$};
\draw [line width=1pt] (0,0)-- (0,-2);
\draw [line width=1pt] (0,0)-- (-3,0);
\draw [line width=1pt] (0.7363201646420022,0.3681600823210011)-- (0.715252084409844,0.41029624278531734);
\draw[->] (2,0) arc (0:30:1.8cm);
\draw (2.,0.7) node[anchor=north west] {$\beta$};
\end{tikzpicture}
         \caption{Case $\frac{\hat{f}_i}{e^{i\beta}}>0$ : $|\hat{g}_i|>|\hat{f}_i|$}
         \label{fig:positive_case}
     \end{subfigure}
     \hfill
     \begin{subfigure}[b]{0.4\textwidth}
         \centering
         
\begin{tikzpicture}
\clip(-3.6593270712778567,-2.6796850995853734) rectangle (4.392471953927827,2.772429935419367);
\draw [->,line width=1pt] (0,0) -- (4,2);
\draw [line width=1pt] (0,0)-- (-2,-1);
\draw [line width=1pt] (-2,-1)-- (-3,-1.5);
\draw [->,line width=1pt] (0,0) -- (4,0);
\draw [->,line width=1pt] (0,0) -- (0,2);
\draw (3.562017345946972,1.8336551611801402) node[anchor=north west] {$e^{i\beta}\mathbb{R}$};
\draw (3.,0.0734524594815899) node[anchor=north west] {$\mathbb{R}$};
\draw (.9,2.1) node[anchor=north west] {$|\hat{f}_i|e^{i\beta}=-\hat{f}_i$};
\draw [line width=1pt] (2.75,1.5)-- (2.78474075768685,1.4210640423044958);
\draw  (-2.13,1.2)node[anchor=north west] {$|\hat{g}_i|e^{i\beta}=(\frac{3}{2}-1)|\hat{f}_i|e^{i\beta}$};
\draw (-0.7,1.9239219663954505) node[anchor=north west] {$i\mathbb{R}$};
\draw [line width=1pt] (0,0)-- (0,-2);
\draw [line width=1pt] (0,0)-- (-3,0);
\draw [line width=1pt] (0.7363201646420022,0.3681600823210011)-- (0.715252084409844,0.41029624278531734);
\draw[->] (2,0) arc (0:30:1.8cm);
\draw (2.,0.7) node[anchor=north west] {$\beta$};
\draw [line width=1pt] (-2.8,-1.4)-- (-2.83,-1.35);
\draw (-3.,-0.5) node[anchor=north west] {$\hat{f}_i$};
\end{tikzpicture}
         \caption{Case $\frac{\hat{f}_i}{e^{i\beta}}<0$ : $|\hat{g}_i|<|\hat{f}_i|$}
         \label{fig:negative_case}
     \end{subfigure}
        \caption{Geometrical representation of sign determination }
        \label{fig:graph_phase}
\end{figure*} Then remains the determination of the different phases of the Fourier coefficients. To do so, a method allowing to distinguish for each coefficients $k\in [0,N-1]$ if  $\frac{\hat{f_k}}{e^{i\beta_k}}$ is positive or negative is needed. Luckily a nice property of the Fourier transform of complex exponential functions can be used to develop a discriminant method. Suppose that there is a known approximation of the amplitudes of the different quantum coefficients obtained by an MLQAE method : $\boldsymbol{\Tilde{a}}=(\Tilde{a}_k)_{k\in[0,N-1]}$ where $\Tilde{a}_k \approx |\hat{f}_k|$. Let's define $\boldsymbol{e^k}=(e^k_p)_{p\in[0,N-1]}\in \C^N$ with $e^k_p=e^{2i\pi kp}$ to ensure that $\hat{\boldsymbol{e}}^k=\boldsymbol{\delta}^k$ where $\delta^k_{p}=1$ if $k=p$ and $\delta^k_p=0$ otherwise. Set $\boldsymbol{g}^k=(g^k_p)_{p\in[0,N-1]} \in \C^N$ such as \begin{equation}
    \boldsymbol{g}^k= \boldsymbol{f}+\frac{3}{2}\Tilde{a}_k \boldsymbol{e}^k e^{\beta_k}.
\end{equation}
This leads to  \begin{equation}
    \hat{\boldsymbol{g}}^k= \hat{\boldsymbol{f}}+\frac{3 }{2}\Tilde{a}_k \boldsymbol{\delta}^k e^{\beta_k}
\end{equation}
In the end there is that \begin{equation}
    \begin{cases}
    \text{if } k \neq j \forall (k,j)\in[0,N-1]^2 \text{ then } \hat{g}^k_j= \hat{f}_j, \\
    \text{otherwise  } \hat{g}^k_k \neq \hat{f}_k.
    \end{cases}
\end{equation}From this two last equation, a sign determination method can be deduced (see Figure \ref{fig:graph_phase}) as we have the property : \begin{equation}
    |\hat{g}^k_k| \geq |\hat{f}_k| \Longleftrightarrow \frac{\hat{f_k}}{e^{\beta_k}} \geq 0.
\end{equation} From this method a QFT based Discrete Fourier Transform  sign determination algorithm can then be written :
\begin{algorithm}[H]
\caption{Sign determination}\label{alg:QFT_based_sign}
\begin{algorithmic}
\Require $\boldsymbol{f},\boldsymbol{M},\boldsymbol{a}$
\For{$k=1,N$}
\State $\boldsymbol{g}^k= \boldsymbol{f}+\frac{3}{2 \sqrt{N}}\Tilde{a}_k \boldsymbol{c}^k e^{\beta_k}.$
\State $\boldsymbol{\mathcal{A}}_k=\boldsymbol{\mathcal{I}}_{\frac{\boldsymbol{g^k}}{\|g^k\|}}\text{QFT}_n\boldsymbol{\mathcal{K}}_k$
\State $\boldsymbol{\mathcal{Q}}_k= \boldsymbol{\mathcal{A}}_k \boldsymbol{\mathcal{S}}_0 \boldsymbol{\mathcal{A}}_k^{-1}\boldsymbol{\mathcal{S}}_{\chi}$
\State $d_k=\sqrt{N} \|g^k\| \mlqae (\boldsymbol{\mathcal{Q}}_k,\boldsymbol{M})$
\If{$d_k \geq a_k$}
    \State $b_k=a_k e^{\beta_k}$
\ElsIf{$d_k < a_k$}
    \State $b_k=-a_k e^{\beta_k}$
\EndIf
\EndFor
\return $\boldsymbol{b}$
\end{algorithmic}
\end{algorithm}
From there the general algorithm of QFT for symmetric input data can be written : 
\begin{algorithm}[H]
\caption{Symmetric QFT }\label{alg:symmetric_QFT}
\begin{algorithmic}
\Require $\boldsymbol{f},\boldsymbol{M}$
\State $\boldsymbol{a}=\text{QFT amplitude}(\boldsymbol{f},\boldsymbol{M})$
\State $\hat{\boldsymbol{f}}=\text{Sign determination}(\boldsymbol{f},\boldsymbol{M},\boldsymbol{a})$
\return $\hat{\boldsymbol{f}}$
\end{algorithmic}
\end{algorithm}
The error of approximation for a given coefficient is following the same asymptotic comportment of the Hadamard test based method (Equation \ref{eqn:asymptotic_error}). This is because the error limiting process is the Most Likelihood Amplitude Estimation also used in the Hadamard test method.  

\subsection{Results}

\subsubsection{Most Likelihood amplitude estimation}

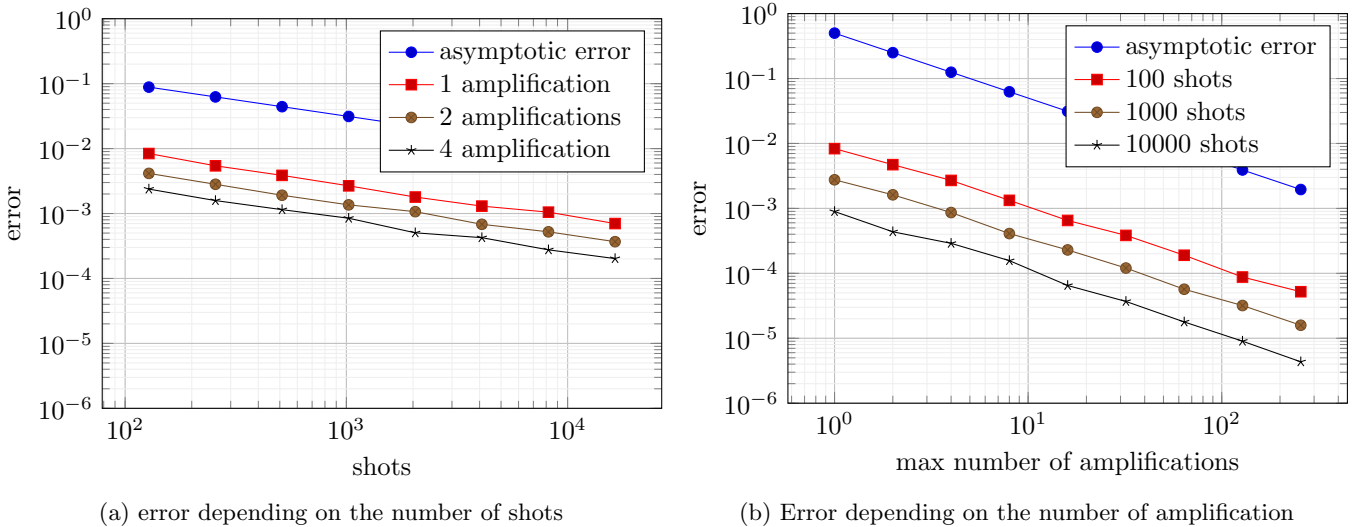
\begin{figure*}
\begin{subfigure}[b]{0.5\textwidth}
    \centering
    \resizebox{\textwidth}{!}{
\begin{tikzpicture}
\begin{axis}[
    legend style={nodes={scale=1, transform shape}},
    xlabel=shots,
    ylabel=error,
    ymode=log,
    xmode=log,
    ymax=1,
    ymin=1e-6,
    grid = both,
    major grid style = {lightgray},
    minor grid style = {lightgray!25},
    width = \textwidth,
    height = 0.75\textwidth,
    legend cell align = {left},
    legend pos = north east,
]
        \addplot table[ x={x},y={th},col sep=comma, mark=* ] {data/simulator/MLQAE/shots_error.dat};
        \addplot table[ x={x},y={amp_1},col sep=comma, mark=* ] {data/simulator/MLQAE/shots_error.dat};
        \addplot table[ x={x},y={amp_2},col sep=comma, mark=* ] {data/simulator/MLQAE/shots_error.dat};
        \addplot table[ x={x},y={amp_3},col sep=comma, mark=* ] {data/simulator/MLQAE/shots_error.dat};

\legend{
asymptotic error,
1 amplification,
2 amplifications,
4 amplification,
}
\end{axis}
\end{tikzpicture}
}
    \caption{error depending on the number of shots}
    \label{fig:MLQAE_shots}
\end{subfigure}
\begin{subfigure}[b]{0.5\textwidth}
    \centering
    \resizebox{\textwidth}{!}{
\begin{tikzpicture}
\begin{axis}[
    legend style={nodes={scale=1, transform shape}},
    xlabel=max number of amplifications,
    ylabel=error,
    ymode=log,
    xmode=log,
    ymax=1,
    ymin=1e-6,
    grid = both,
    major grid style = {lightgray},
    minor grid style = {lightgray!25},
    width = \textwidth,
    height = 0.75\textwidth,
    legend cell align = {left},
    legend pos = north east,
]
        \addplot table[ x={x},y={th},col sep=comma, mark=* ] {data/simulator/MLQAE/amp_error.dat};
        \addplot table[ x={x},y={ml1},col sep=comma, mark=* ] {data/simulator/MLQAE/amp_error.dat};
        \addplot table[ x={x},y={ml2},col sep=comma, mark=* ] {data/simulator/MLQAE/amp_error.dat};
        \addplot table[ x={x},y={ml3},col sep=comma, mark=* ] {data/simulator/MLQAE/amp_error.dat};

\legend{
asymptotic error,
100 shots,
1000 shots,
10000 shots,
}
\end{axis}
\end{tikzpicture}
}
    \caption{Error depending on the number of amplification }
    \label{fig:MLQAE_amp}
\end{subfigure}
    \caption{MLQAE error analysis}
    \label{fig:MLQA_analysis}
\end{figure*}

As a first test, the MLQAE method has been performed on a very simple input initialized on 1 qubit in order to refine the asymptotic comportment of the error curves. To do so, the input vector chosen is  $f=\begin{pmatrix} \frac{1}{\sqrt{2}}\\ \frac{1}{\sqrt{2}}
\end{pmatrix}$. This state is chosen because it maximizes the function $s(x)=\sqrt{|x|^2(1-|x|^2)}$. It is known from \cite{suzuki2020amplitude} that to estimate the amplitude of a quantum coefficient $x_k$, MLQAE error has a scaling factor of $s(x_k)$. This state can be initialized using a $\boldsymbol{H}$ gate on a $\ket{0}{}$ state, so we can deduce $\boldsymbol{\mathcal{A}}_k=\boldsymbol{H}\boldsymbol{\mathcal{K}}_k$. The test has been performed in two different ways, first by using fixed number of amplification and by making the number of repetitions of the experiment variate. Second by using fixed number of repetitions and by making the number of amplification variate. By this, it will be possible to verify and maybe refine the asymptotic comportment of the error of MLQAE method given by equation \ref{eqn:asymptotic_error} to ensure limited stochastic variations of the results due to the character of the MLQAE method, the experiments have been repeated 50 times. From the results it can be observed that the asymptotic comportment is in accordance with the theory. The shots varying experiment (Figure \ref{fig:MLQAE_shots}) shows that for all the number of amplification tested the decrease of the error is done in $O(\frac{1}{\sqrt{N}})$ where $N$ is the number of repetition. More precisely, a linear regression on logarithmic scale shows a slope of $0.5 \pm 0.01$ with a correlation factor of $0.99$ and a $p$ value of $10e-12$. In the case of number of amplifications varying similar results can be observed, here for all the number of repetitions tested the error decreases in $O(\frac{1}{m})$ where $m$ is the maximum number of amplifications. More precisely, a linear regression on logarithmic scale shows a slope of $0.95 \pm 0.06$ with a correlation factor $0.99$ and a $p$ value $10e-12$. Even if there is a little gap with the asymptotic slope the results is still in the deviation interval which makes the two results are in relative accordance with the theoretical asymptotic comportment of MLQAE method. nevertheless it could be refined in order to make possible a way of having an \emph{a priori} controlled error by choosing the number of amplifications in order to have a targeted error. This could be done by extending Equation \ref{eqn:asymptotic_error}, as there exist $C \in R$ and $r$ a function from $\mathbb{N}^2$ to $\R$ negligible \begin{equation}
    e(m_M,S)=\frac{C}{m_M\sqrt{S}}+r(M,S).
\end{equation} An expression of $m_M$ depending on a targeted error can be found : \begin{equation}
    m_M(S,e)=\mathcal{O}(\frac{1}{e\sqrt{S}}). 
    \label{eqn:amplifcation_number}
\end{equation} Nevertheless this expression is far from being precise, in fact if we use it do deduce the number of amplification needed to measure our input vector coefficients with an accuracy of $10e-4$ and 1000 shots, the number of amplification necessary would be around $300$ which is far from what can be observed on Figure \ref{fig:MLQAE_amp}. This gap can nevertheless be reduced by estimating the difference between MLQAE error curves and the theoretical asymptotic error curve which will give an estimation for the constant $C$. This can be done by taking the intersection of the error curves $\text{int}$ with the error axis in logarithmic scale in the Figure \ref{fig:MLQAE_shots} and multiplying by the number of amplification : $C=\text{int} \times m$. On Figure \ref{fig:MLQAE_amp} it can also be done by taking the intersection between the error curves and the error axis and multiplying by the square root of the number of repetitions : $C=\text{int}\sqrt{S}$. From the experiments of Figure\ref{fig:MLQA_analysis} an approximation of $C$ is calculated, the results are summarized in table. A refined asymptotic comportment of m depending on the error can now be defined by \begin{equation}
    m_M(S,e)=O(\frac{C}{e\sqrt{S}}),
\end{equation} where the value of $C$ can be chosen among the one measured. Using this formula to determine the number of shots necessary to have an error of approximation of $10e-4$ gives now a way more accurate result around $32$ comparable to what can be observed on Figure \ref{fig:MLQAE_amp}.

\begin{table}[h!]
\centering
\begin{tabular}{ |c|c|c| } 

\hline
 method & $C$ & standard deviation  \\
\hline
 shots & 0.088 & 0.009  \\
amplifications  & 0.094 & 0.003  \\
\hline
\end{tabular}
\caption{Table of values of the constant $C$}
\label{table:2}
\end{table}

\subsubsection{most likelihood amplitude estimation results}

In this section, the same tests have been performed on real devices to see if the same comportment was observed. Using the same method than in the last section, an input vector  $f=\begin{pmatrix} \frac{1}{\sqrt{2}}\\ \frac{1}{\sqrt{2}}
\end{pmatrix}$ is used, as it only needs an $\boldsymbol{H}$ gate as initialization. The results observed on IBMQ Ehningen are shown on Figure \ref{fig:MLQA_analysis_lima}. The analysis on the number of shots is showing that the number of shots doesn't have a clear impact on the results which is counter intuitive. Fortunately, it seems that the number amplifications allows to reach more accuracy. With the analysis on the number of amplification it can be observed that the curve is not exactly linear, more precisely the p value of the linear regression is greater than 0.05 which reject this hypothesis. This comportment can be explained by several facts. First measurement error (predominant here) around $1e-2$ slows down MLQAE convergence responsible of the bump curves that can be observed. To solve this problem a measurement error mitigation is used to correct this type of noise. Results can be observed on Figure \ref{fig:MLQA_analysis_lima_mitigated}, it can be be noticed that the bump observed on Figure \ref{fig:MLQA_analysis_lima} is less present. Especially for the 128 shots curves it can be seen that the comportment became quasi linear with a decreasing slope of coefficient $-0.87$. This slight difference in the asymptotic comportment can be explained by the other types of noise such as decoherence error or gates errors. Unfortunately this errors types are complicated to mitigate.

\begin{figure*}
\begin{subfigure}[b]{0.5\textwidth}
    \centering
    \resizebox{\textwidth}{!}{
\begin{tikzpicture}
\begin{axis}[
    legend style={nodes={scale=1, transform shape}},
    xlabel=shots,
    ylabel=error,
    ymode=log,
    xmode=log,
    ymax=1,
    ymin=1e-6,
    grid = both,
    major grid style = {lightgray},
    minor grid style = {lightgray!25},
    width = \textwidth,
    height = 0.75\textwidth,
    legend cell align = {left},
    legend pos = north east,
]
        \addplot table[ x={x},y={th},col sep=comma, mark=* ] {data/simulator/MLQAE/shots_error.dat};
        \addplot table[ x={x},y={amp_1},col sep=comma, mark=* ] {data/real_devices/ehningen/shots_error.dat};
        \addplot table[ x={x},y={amp_8},col sep=comma, mark=* ] {data/real_devices/ehningen/shots_error.dat};
        \addplot table[ x={x},y={amp_9},col sep=comma, mark=* ] {data/real_devices/ehningen/shots_error.dat};
        \addplot table[ x={x},y={amp_10},col sep=comma, mark=* ] {data/real_devices/ehningen/shots_error.dat};

\legend{
asymptotic error,
1 amplification,
256 amplifications,
512 amplifications,
1024 amplifications,
}
\end{axis}
\end{tikzpicture}
}
    \caption{Error depending on the number of shots}
    \label{fig:MLQAE_shots}
\end{subfigure}
\begin{subfigure}[b]{0.5\textwidth}
    \centering
    \resizebox{\textwidth}{!}{
\begin{tikzpicture}
\begin{axis}[
    legend style={nodes={scale=1, transform shape}},  
    xlabel=max number of amplifications,
    ylabel=error,
    ymode=log,
    xmode=log,
    grid = both,
    ymax=1,
    ymin=1e-6,
    major grid style = {lightgray},
    minor grid style = {lightgray!25},
    width = \textwidth,
    height = 0.75\textwidth,
    legend cell align = {left},
    legend pos = north east,
]
        \addplot table[ x={x},y={th},col sep=comma, mark=* ] {data/simulator/MLQAE/amp_error.dat};
        \addplot table[ x={x},y={amp1},col sep=comma, mark=* ] {data/real_devices/ehningen/amp_error.dat};
        \addplot table[ x={x},y={amp2},col sep=comma, mark=* ] {data/real_devices/ehningen/amp_error.dat};
        \addplot table[ x={x},y={amp3},col sep=comma, mark=* ] {data/real_devices/ehningen/amp_error.dat};

\legend{
asymptotic error,
100 shots,
1000 shots,
10000 shots,
}
\end{axis}
\end{tikzpicture}
}
    \caption{Error depending on the number of amplification }
    \label{fig:MLQAE_amp}
\end{subfigure}
    \caption{MLQAE error analysis on IBMQ Ehningen}
    \label{fig:MLQA_analysis_lima}
\end{figure*}

\begin{figure*}
\begin{subfigure}[b]{0.5\textwidth}
    \centering
    \resizebox{\textwidth}{!}{
\begin{tikzpicture}
\begin{axis}[
    legend style={nodes={scale=1, transform shape}},
    xlabel=shots,
    ylabel=error,
    ymode=log,
    xmode=log,
    ymax=1,
    ymin=1e-6,
    grid = both,
    major grid style = {lightgray},
    minor grid style = {lightgray!25},
    width = \textwidth,
    height = 0.75\textwidth,
    legend cell align = {left},
    legend pos = north east,
]
        \addplot table[ x={x},y={th},col sep=comma, mark=* ] {data/simulator/MLQAE/shots_error.dat};
        \addplot table[ x={x},y={amp_1},col sep=comma, mark=* ] {data/real_devices/ehningen/mitigated/shots_error.dat};
        \addplot table[ x={x},y={amp_8},col sep=comma, mark=* ] {data/real_devices/ehningen/mitigated/shots_error.dat};
        \addplot table[ x={x},y={amp_9},col sep=comma, mark=* ] {data/real_devices/ehningen/mitigated/shots_error.dat};
        \addplot table[ x={x},y={amp_10},col sep=comma, mark=* ] {data/real_devices/ehningen/mitigated/shots_error.dat};

\legend{
asymptotic error,
1 amplification,
256 amplifications,
512 amplifications,
1024 amplifications,
}
\end{axis}
\end{tikzpicture}
}
    \caption{Error depending on the number of shots}
    \label{fig:MLQAE_shots}
\end{subfigure}
\begin{subfigure}[b]{0.5\textwidth}
    \centering
    \resizebox{\textwidth}{!}{
\begin{tikzpicture}
\begin{axis}[
    legend style={nodes={scale=1, transform shape}},  
    xlabel=max number of amplifications,
    ylabel=error,
    ymode=log,
    xmode=log,
    grid = both,
    ymax=1,
    ymin=1e-6,
    major grid style = {lightgray},
    minor grid style = {lightgray!25},
    width = \textwidth,
    height = 0.75\textwidth,
    legend cell align = {left},
    legend pos = north east,
]
        \addplot table[ x={x},y={th},col sep=comma, mark=* ] {data/real_devices/ehningen/mitigated/amp_error.dat};
        \addplot table[ x={x},y={amp_1},col sep=comma, mark=* ] {data/real_devices/ehningen/mitigated/amp_error.dat};
        \addplot table[ x={x},y={amp_5},col sep=comma, mark=* ] {data/real_devices/ehningen/mitigated/amp_error.dat};
        \addplot table[ x={x},y={amp_7},col sep=comma, mark=* ] {data/real_devices/ehningen/mitigated/amp_error.dat};

\legend{
asymptotic error,
128 shots,
1024 shots,
8192 shots,
}
\end{axis}
\end{tikzpicture}
}
    \caption{Error depending on the number of amplification }
    \label{fig:MLQAE_amp}
\end{subfigure}
    \caption{MLQAE error analysis on IBMQ Ehningen}
    \label{fig:MLQA_analysis_lima_mitigated}
\end{figure*}

\subsubsection{QFT results}

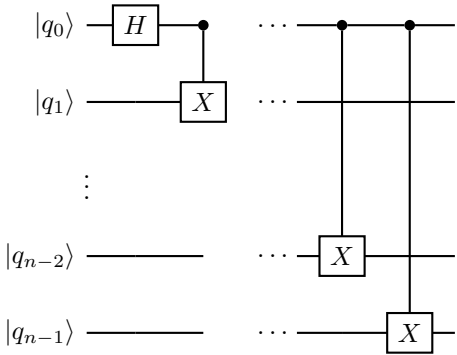
\begin{figure}[H]
\centering
\begin{tikzpicture}
\node[scale=1] {
\begin{quantikz}[column sep=0.3cm]
\lstick{$\ket{q_0}$} & \gate{H} &   \ctrl{1} & \ \dots \ & \ctrl{3} & \ctrl{4} &  \qw \\
\lstick{$\ket{q_1}$} & \qw & \gate{X} & \ \dots \ & \qw & \qw & \qw \\
\vdots &&&&&&&&&&&&&&&&  \\ 
\lstick{$\ket{q_{n-2}}$}& \qw & \qw & \ \dots \  & \gate{X} & \qw & \qw \\
\lstick{$\ket{q_{n-1}}$} &\qw & \qw & \ \dots \ & \qw& \gate{X} & \qw 
\end{quantikz}
};
\end{tikzpicture}
    \caption{Bell state Initialization circuit}
    \label{fig:bellstateini}
\end{figure}

\begin{figure*}
\begin{subfigure}[b]{0.5\textwidth}
    \centering
    \resizebox{0.8\textwidth}{!}{
    \begin{tikzpicture}
\begin{axis}[
    legend style={nodes={scale=1, transform shape}},
    xlabel=$\xi$,
    ylabel=$\frac{|\hat{f}|}{\|\hat{f}\|}$,
    xmin = 0, xmax = 3,
    ymin = 0, ymax = 4,
    xtick distance = 1,
    ytick distance = 1,
    grid = both,
    minor tick num = 1,
    major grid style = {lightgray},
    minor grid style = {lightgray!25},
    width = \textwidth,
    height = 0.75\textwidth,
    legend cell align = {left},
    legend pos = north west
]
 
        \addplot table[x={x},y expr=4*\thisrow{q0},col sep=comma, mark= *] {\opbt};
        \addplot table[x={x},y expr=4*\thisrow{q1},col sep=comma, mark= *] {\opbt};
        \addplot table[x={x},y expr=4*\thisrow{q2},col sep=comma, mark= *] {\opbt};
        \addplot table[x={x},y expr=4*\thisrow{q3},col sep=comma, mark= *] {\opbt};

        \addplot table[color=red, x={x},y expr=4*\thisrow{rv},col sep=comma, mark= ] {\opbt};
 
\legend{
    QFT amplification order : 0, 
    QFT amplification order : 1, 
    QFT amplification order : 2, 
    QFT amplification order : 4, 
    FFT
}
 
\end{axis}
 
\end{tikzpicture}
}
    \caption{2 qubits 1D QFT}
    \label{fig:my_label}
\end{subfigure}
\begin{subfigure}[b]{0.5\textwidth}
    \centering
    \resizebox{0.8\textwidth}{!}{
    \begin{tikzpicture}
\begin{axis}[
    legend style={nodes={scale=1, transform shape}},
    xlabel=$\xi$,
    ylabel=$\frac{|\hat{f}|}{\|\hat{f}\|}$,
    xmin = 0, xmax = 7,
    ymin = 0, ymax = 4,
    xtick distance = 1,
    ytick distance = 1,
    grid = both,
    minor tick num = 1,
    major grid style = {lightgray},
    minor grid style = {lightgray!25},
    width = \textwidth,
    height = 0.75\textwidth,
    legend cell align = {left},
    legend pos = north west
]
 
        \addplot table[x={x},y expr=8*\thisrow{q0},col sep=comma, mark= *] {\opbtr};
        \addplot table[x={x},y expr=8*\thisrow{q1},col sep=comma, mark= *] {\opbtr};
        \addplot table[x={x},y expr=8*\thisrow{q2},col sep=comma, mark= *] {\opbtr};
        \addplot table[x={x},y expr=8*\thisrow{q3},col sep=comma, mark= *] {\opbtr};
        \addplot table[color=red, x={x},y expr=8*\thisrow{rv},col sep=comma, mark= ] {\opbtr};
 
\legend{
    QFT amplification order : 0, 
    QFT amplification order : 1, 
    QFT amplification order : 2, 
    QFT amplification order : 4, 
    FFT
}
 
\end{axis}
 
\end{tikzpicture}
}
    \caption{3 qubits 1D QFT}
    \label{fig:my_label}
\end{subfigure}

\begin{subfigure}[b]{0.5\textwidth}
    \centering
    \resizebox{0.8\textwidth}{!}{
    \begin{tikzpicture}
\begin{axis}[
    legend style={nodes={scale=1, transform shape}},
    xlabel=$\xi$,
    ylabel=$\frac{|\hat{f}|}{\|\hat{f}\|}$,
    xmin = 0, xmax = 15,
    ymin = 0, ymax = 4,
    xtick distance = 1,
    ytick distance = 1,
    grid = both,
    minor tick num = 1,
    major grid style = {lightgray},
    minor grid style = {lightgray!25},
    width = \textwidth,
    height = 0.75\textwidth,
    legend cell align = {left},
    legend pos = north west
]
 
        \addplot table[x={x},y expr=16*\thisrow{q0},col sep=comma, mark= *] {\opbf};
        \addplot table[x={x},y expr=16*\thisrow{q1},col sep=comma, mark= *] {\opbf};
        \addplot table[x={x},y expr=16*\thisrow{q2},col sep=comma, mark= *] {\opbf};
        \addplot table[x={x},y expr=16*\thisrow{q3},col sep=comma, mark= *] {\opbf};
        \addplot table[color=red, x={x},y expr=16*\thisrow{rv},col sep=comma, mark= ] {\opbf};
 
\legend{
    QFT amplification order : 0, 
    QFT amplification order : 1, 
    QFT amplification order : 2, 
    QFT amplification order : 4, 
    FFT
}
 
\end{axis}
 
\end{tikzpicture}
}
    \caption{4 qubits 1D QFT}
    \label{fig:my_label}
\end{subfigure}
\begin{subfigure}[b]{0.5\textwidth}
    \centering
    \resizebox{0.8\textwidth}{!}{
    \begin{tikzpicture}
\begin{axis}[
    legend style={nodes={scale=1, transform shape}},
    xlabel=$\xi$,
    ylabel=$\frac{|\hat{f}|}{\|\hat{f}\|}$,
    xmin = 0, xmax = 31,
    ymin = 0, ymax = 4,
    xtick distance = 2,
    ytick distance = 1,
    grid = both,
    minor tick num = 1,
    major grid style = {lightgray},
    minor grid style = {lightgray!25},
    width = \textwidth,
    height = 0.75\textwidth,
    legend cell align = {left},
    legend pos = north west
]
 
        \addplot table[x={x},y expr=32*\thisrow{q0},col sep=comma, mark= *] {\opbfi};
        \addplot table[x={x},y expr=32*\thisrow{q1},col sep=comma, mark= *] {\opbfi};
        \addplot table[x={x},y expr=32*\thisrow{q2},col sep=comma, mark= *] {\opbfi};
        \addplot table[x={x},y expr=32*\thisrow{q3},col sep=comma, mark= *] {\opbfi};
        \addplot table[color=red, x={x},y expr=32*\thisrow{rv},col sep=comma, mark= ] {\opbfi};
 
\legend{
    QFT amplification order : 0, 
    QFT amplification order : 1, 
    QFT amplification order : 2, 
    QFT amplification order : 4, 
    FFT
}
 
\end{axis}
 
\end{tikzpicture}
}
    \caption{5 qubits 1D QFT}
    \label{fig:my_label}
\end{subfigure}

    \caption{Bell state Input Quantum Fourier transform using Hadamard test}
    \label{fig:HADAMARD1D}
\end{figure*}
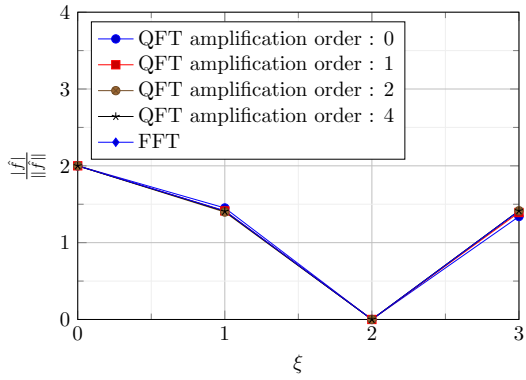
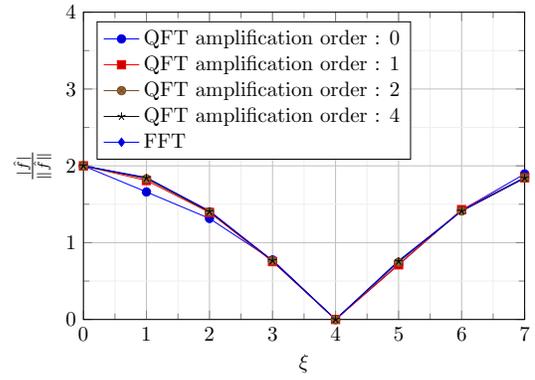
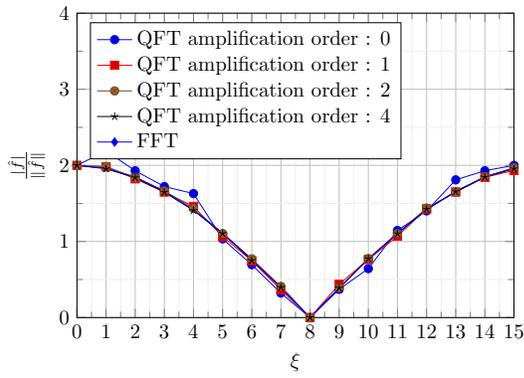
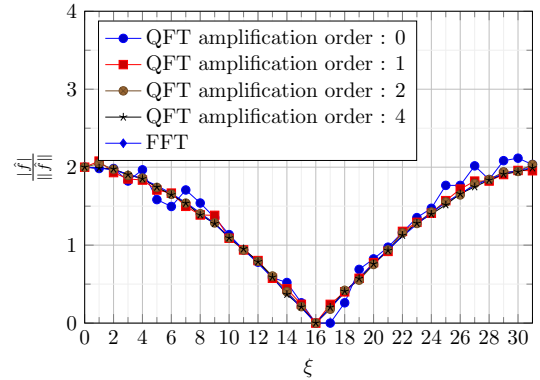

As a first step we applied the Quantum Fourier transform in dimension 1 in order to compare the two methods of Fourier coefficients estimation. As test cases, for a given number of qubits $n$, the input vector $f \in \R^{2^n}$ will be of the form $$f=\begin{pmatrix} \frac{1}{\sqrt{2}}\\ 0\\ \vdots\\ 0\\ \frac{1}{\sqrt{2}}
\end{pmatrix}.$$ This input state is called a Bell state. It has several advantages, first it's symmetric, so both methods will give exact complex Fourier coefficients, moreover it's an input quantum state with a pretty simple initialization circuit (as shown in Figure \ref{fig:bellstateini}) so we can expect shorter Initialization circuits compared to the ones calculated with the method of Mottonen \cite{mottonen2004transformation}. To analyze the different results, the $L^2$ error is used, it is defined by \begin{equation*}
    \|f_{\text{exact}}-f_{\text{estimated}}\|_{L^2} =\sqrt{ \sum^{2^n-1}_{k=0}(f_{\text{exact}_k}-f_{\text{estimated}_k})^2}.
\end{equation*} As first results some tests have been performed on Qiskit QASM simulator for different number of qubits and different number of amplifications (Figure \ref{fig:HADAMARD1D} and Figure \ref{fig:SIGN1D}). The number of shots is fixed at 1000 as it has been seen on real devices that after 1000, the number of shots had no significant effect on the error. This simulator is acting as a "perfect" quantum device in the sense that it doesn't introduce any noise which allows to test if the algorithms are working properly in a first place. As a first conclusion from this results, it can be seen that the two algorithms are acting similarly with the number of amplifications. In fact it can be seen that they follow a law in $\frac{1}{m_M}$ which is coherent with what been seen in the method part. It can also be seen that the two methods have comparable error comportment (Figure \ref{fig:L21D}). Then, the advantage of the sign determination method comes from the fact that it has a shorter number of gates (Figure \ref{fig:circuit_length}) which may lead to a better robustness to circuit length related noise and decoherence noise on real devices. Nonetheless, it can be seen that the error is also increasing with the number of qubits, this come from the fact that the estimation is done on the normalize Fourier transform. In fact increasing the size of the domain comes in general with an increasing norm of the Fourier transform which makes the error rescaled in consequence. 
\newline

However, when it comes to real devices some tests have been performed on different different NISQ devices, the result showed are the one given by the German Quantum computer IBMQ Ehningen. As a first result with 2 qubits, it can be seen that even without amplification Hadamard test algorithms are giving completely noisy results. This can be explain by several facts. First, the fact that the controlled operation in the Hadamard test circuits is composed with a lot of Controlled NOT gates which are susceptible to noise (probability of success around $1e-2$ compared to $1e-4$ for any other gate). Second, the Hadamard test circuit needs one more qubit to work than the sign estimation method and as all the qubits are not connected together some swap gates have to be used during the calculation, this swap gates are also composed of 3 CNOT gates which leeds to an even higher number of CNOT and then an higher gate related noise. A third problem is also the decoherence error due to the interaction of the system with the environment which collapse the wave function. It is hard to estimate in practice but it's modeled by a probability of a bit flip and a phase flip increasing with the time of the calculation and so the length of the circuit. In the other hand, the results obtained for the sign determination methods are not really better, even if we see that the results are not completely noisy, they are also impossible to use in practice as the noise is still very present. Even more, results with a higher number of qubits show that even without amplification the noise predominates. And even if in the sign determination algorithm CNOT gates are less used, there is still a lot due to the QFT and some other controlled operation (including the initialization) which could explain this high noise compared to some results on real devices found in the literature for this type of amplitude estimation.  


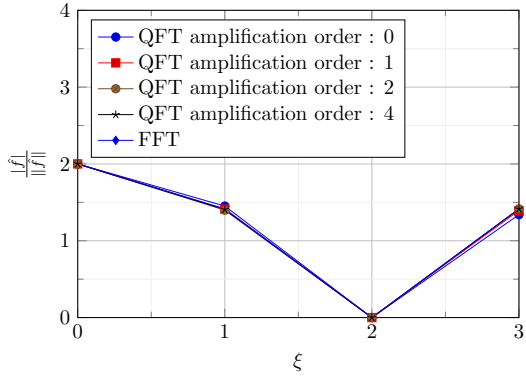
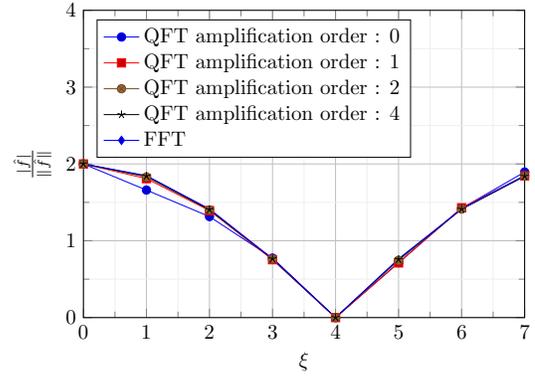
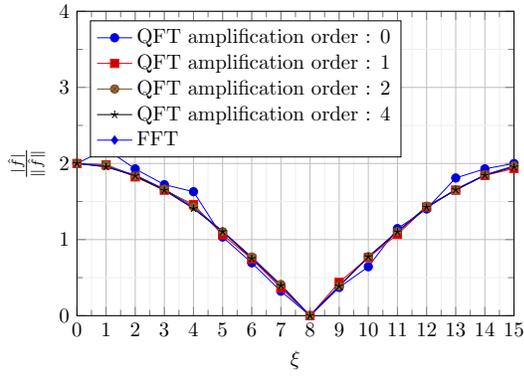
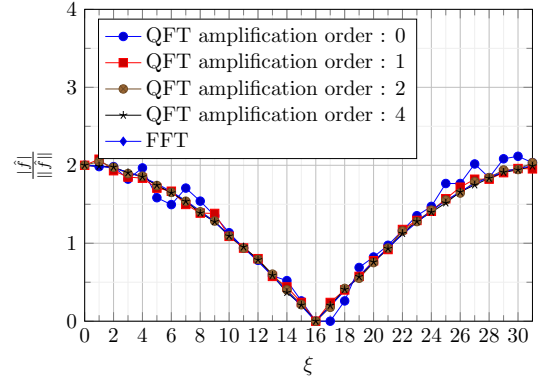
\begin{figure*}
\begin{subfigure}[b]{0.5\textwidth}
    \centering
        \resizebox{0.8\textwidth}{!}{
    \begin{tikzpicture}
\begin{axis}[
    legend style={nodes={scale=1, transform shape}},
    xlabel=$\xi$,
    ylabel=$\frac{|\hat{f}|}{\|\hat{f}\|}$,
    xmin = 0, xmax = 3,
    ymin = 0, ymax = 4,
    xtick distance = 1,
    ytick distance = 1,
    grid = both,
    minor tick num = 1,
    major grid style = {lightgray},
    minor grid style = {lightgray!25},
    width = \textwidth,
    height = 0.75\textwidth,
    legend cell align = {left},
    legend pos = north west
]
 
        \addplot table[x={x},y expr=4*\thisrow{q0},col sep=comma, mark= *] {\opbt};
        \addplot table[x={x},y expr=4*\thisrow{q1},col sep=comma, mark= *] {\opbt};
        \addplot table[x={x},y expr=4*\thisrow{q2},col sep=comma, mark= *] {\opbt};
        \addplot table[x={x},y expr=4*\thisrow{q3},col sep=comma, mark= *] {\opbt};

        \addplot table[color=red, x={x},y expr=4*\thisrow{rv},col sep=comma, mark= ] {\opbt};
 
\legend{
    QFT amplification order : 0, 
    QFT amplification order : 1, 
    QFT amplification order : 2, 
    QFT amplification order : 4, 
    FFT
}
 
\end{axis}
 
\end{tikzpicture}
}
    \caption{2 qubits 1D QFT}
    \label{fig:my_label}
\end{subfigure}
\begin{subfigure}[b]{0.5\textwidth}
    \centering
        \resizebox{0.8\textwidth}{!}{
    \begin{tikzpicture}
\begin{axis}[
    legend style={nodes={scale=1, transform shape}},
    xlabel=$\xi$,
    ylabel=$\frac{|\hat{f}|}{\|\hat{f}\|}$,
    xmin = 0, xmax = 7,
    ymin = 0, ymax = 4,
    xtick distance = 1,
    ytick distance = 1,
    grid = both,
    minor tick num = 1,
    major grid style = {lightgray},
    minor grid style = {lightgray!25},
    width = \textwidth,
    height = 0.75\textwidth,
    legend cell align = {left},
    legend pos = north west
]
 
        \addplot table[x={x},y expr=8*\thisrow{q0},col sep=comma, mark= *] {\opbtr};
        \addplot table[x={x},y expr=8*\thisrow{q1},col sep=comma, mark= *] {\opbtr};
        \addplot table[x={x},y expr=8*\thisrow{q2},col sep=comma, mark= *] {\opbtr};
        \addplot table[x={x},y expr=8*\thisrow{q3},col sep=comma, mark= *] {\opbtr};
        \addplot table[color=red, x={x},y expr=8*\thisrow{rv},col sep=comma, mark= ] {\opbtr};
 
\legend{
    QFT amplification order : 0, 
    QFT amplification order : 1, 
    QFT amplification order : 2, 
    QFT amplification order : 4, 
    FFT
}
 
\end{axis}
 
\end{tikzpicture}
}
    \caption{3 qubits 1D QFT}
    \label{fig:my_label}
\end{subfigure}

\begin{subfigure}[b]{0.5\textwidth}
    \centering
        \resizebox{0.8\textwidth}{!}{
    \begin{tikzpicture}
\begin{axis}[
    legend style={nodes={scale=1, transform shape}},
    xlabel=$\xi$,
    ylabel=$\frac{|\hat{f}|}{\|\hat{f}\|}$,
    xmin = 0, xmax = 15,
    ymin = 0, ymax = 4,
    xtick distance = 1,
    ytick distance = 1,
    grid = both,
    minor tick num = 1,
    major grid style = {lightgray},
    minor grid style = {lightgray!25},
    width = \textwidth,
    height = 0.75\textwidth,
    legend cell align = {left},
    legend pos = north west
]
 
        \addplot table[x={x},y expr=16*\thisrow{q0},col sep=comma, mark= *] {\opbf};
        \addplot table[x={x},y expr=16*\thisrow{q1},col sep=comma, mark= *] {\opbf};
        \addplot table[x={x},y expr=16*\thisrow{q2},col sep=comma, mark= *] {\opbf};
        \addplot table[x={x},y expr=16*\thisrow{q3},col sep=comma, mark= *] {\opbf};
        \addplot table[color=red, x={x},y expr=16*\thisrow{rv},col sep=comma, mark= ] {\opbf};
 
\legend{
    QFT amplification order : 0, 
    QFT amplification order : 1, 
    QFT amplification order : 2, 
    QFT amplification order : 4, 
    FFT
}
 
\end{axis}
 
\end{tikzpicture}
}
    \caption{4 qubits 1D QFT}
    \label{fig:my_label}
\end{subfigure}
\begin{subfigure}[b]{0.5\textwidth}
    \centering
        \resizebox{0.8\textwidth}{!}{
    \begin{tikzpicture}
\begin{axis}[
    legend style={nodes={scale=1, transform shape}},
    xlabel=$\xi$,
    ylabel=$\frac{|\hat{f}|}{\|\hat{f}\|}$,
    xmin = 0, xmax = 31,
    ymin = 0, ymax = 4,
    xtick distance = 2,
    ytick distance = 1,
    grid = both,
    minor tick num = 1,
    major grid style = {lightgray},
    minor grid style = {lightgray!25},
    width = \textwidth,
    height = 0.75\textwidth,
    legend cell align = {left},
    legend pos = north west
]
 
        \addplot table[x={x},y expr=32*\thisrow{q0},col sep=comma, mark= *] {\opbfi};
        \addplot table[x={x},y expr=32*\thisrow{q1},col sep=comma, mark= *] {\opbfi};
        \addplot table[x={x},y expr=32*\thisrow{q2},col sep=comma, mark= *] {\opbfi};
        \addplot table[x={x},y expr=32*\thisrow{q3},col sep=comma, mark= *] {\opbfi};
        \addplot table[color=red, x={x},y expr=32*\thisrow{rv},col sep=comma, mark= ] {\opbfi};
 
\legend{
    QFT amplification order : 0, 
    QFT amplification order : 1, 
    QFT amplification order : 2, 
    QFT amplification order : 4, 
    FFT
}
 
\end{axis}
 
\end{tikzpicture}
}
    \caption{5 qubits 1D QFT}
    \label{fig:my_label}
\end{subfigure}

    \caption{Bell state Input Quantum Fourier transform using sign determination methods}
    \label{fig:SIGN1D}
\end{figure*}
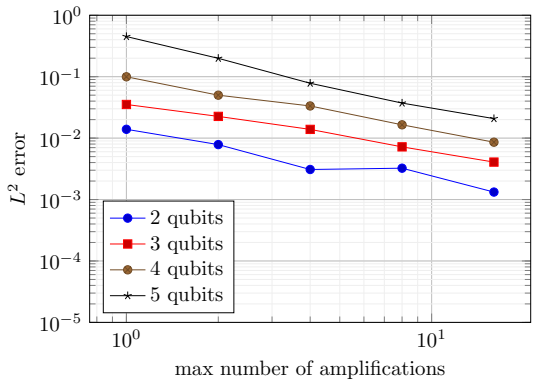
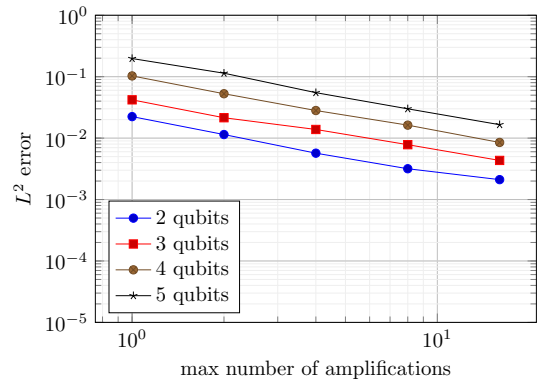
\begin{figure*}
\begin{subfigure}[b]{0.5\textwidth}
    \centering
\resizebox{0.8\textwidth}{!}{
\begin{tikzpicture}
\begin{axis}[
    legend style={nodes={scale=1, transform shape}},
    ymode=log,
    ymax=1,ymin=1e-5,
    xmode=log,
    xlabel=max number of amplifications,
    ylabel=$L^2$ error,
    grid = both,
    minor tick num = 1,
    major grid style = {lightgray},
    minor grid style = {lightgray!25},
    width = \textwidth,
    height = 0.75\textwidth,
    legend cell align = {left},
    legend pos = south west,
]
        \addplot table[ x={x},y expr=\thisrow{q2}^(1/2),col sep=comma, mark=* ] {data/simulator/1D/phase/bell_state/L2/L2.dat};
        \addplot table[ x={x},y expr=\thisrow{q3}^(1/2),col sep=comma, mark=* ] {data/simulator/1D/phase/bell_state/L2/L2.dat};
        \addplot table[ x={x},y expr=\thisrow{q4}^(1/2),col sep=comma, mark=* ] {data/simulator/1D/phase/bell_state/L2/L2.dat};
        \addplot table[ x={x},y expr=\thisrow{q5}^(1/2),col sep=comma, mark=* ] {data/simulator/1D/phase/bell_state/L2/L2.dat};

\legend{2 qubits,
3 qubits,
4 qubits,
5 qubits,}
\end{axis}
\end{tikzpicture}
}
    \caption{Symmetric input method}
    \label{fig:L21DHadamard}
\end{subfigure}
\begin{subfigure}[b]{0.5\textwidth}
    \centering
    \resizebox{0.8\textwidth}{!}{
\begin{tikzpicture}
\begin{axis}[
    legend style={nodes={scale=1, transform shape}},    
    ymode=log,
    ymax=1,ymin=1e-5,
    xmode=log,
    xlabel=max number of amplifications,
    ylabel=$L^2$ error,
    grid = both,
    minor tick num = 1,
    major grid style = {lightgray},
    minor grid style = {lightgray!25},
    width = \textwidth,
    height = 0.75\textwidth,
    legend cell align = {left},
    legend pos = south west,
]
        \addplot table[ x={x},y expr=\thisrow{q2}^(1/2),col sep=comma, mark=* ] {data/simulator/1D/Hadamard/bell_state/L2/L2.dat};
        \addplot table[ x={x},y expr=\thisrow{q3}^(1/2),col sep=comma, mark=* ] {data/simulator/1D/Hadamard/bell_state/L2/L2.dat};
        \addplot table[ x={x},y expr=\thisrow{q4}^(1/2),col sep=comma, mark=* ] {data/simulator/1D/Hadamard/bell_state/L2/L2.dat};
        \addplot table[ x={x},y expr=\thisrow{q5}^(1/2),col sep=comma, mark=* ] {data/simulator/1D/Hadamard/bell_state/L2/L2.dat};
\legend{2 qubits,
3 qubits,
4 qubits,
5 qubits,}
\end{axis}
\end{tikzpicture}
}
    \caption{Hadamard test}
    \label{fig:L21DSIGN}
\end{subfigure}
    \caption{$L_2$ error depending on the number of amplification}
    \label{fig:L21D}
\end{figure*}
\begin{figure}
    \centering
    \resizebox{0.5\textwidth}{!}{
\begin{tikzpicture}
\begin{axis}[
    legend style={nodes={scale=1.2, transform shape}},
    xlabel=max number of amplifications,
    ylabel=circuit length ,
    xmode=log,
    ymode=log,
    grid = both,
    minor tick num = 1,
    major grid style = {lightgray},
    minor grid style = {lightgray!25},
    width = \textwidth,
    height = 0.75\textwidth,
    legend cell align = {left},
    legend pos = north west,
    ]
        \addplot+[blue,mark=square*] table[ x={x},y={q2},col sep=comma, mark=* ] {data/simulator/1D/phase/bell_state/length_circuits/l2qb.dat};
        \addplot+[red,mark=square*] table[ x={x},y={q3},col sep=comma, mark=* ] {data/simulator/1D/phase/bell_state/length_circuits/l2qb.dat};
        \addplot+[brown,mark=square*] table[ x={x},y={q4},col sep=comma, mark=* ] {data/simulator/1D/phase/bell_state/length_circuits/l2qb.dat};
        \addplot+[black,mark=square*] table[ x={x},y={q5},col sep=comma, mark=* ] {data/simulator/1D/phase/bell_state/length_circuits/l2qb.dat};
        \addplot+[blue,dashed,mark=*] table[ x={x},y={q2},col sep=comma, mark=* ] {data/simulator/1D/Hadamard/bell_state/length_circuits/l2qb.dat};
        \addplot+[red,dashed,mark=*] table[ x={x},y={q3},col sep=comma, mark=* ] {data/simulator/1D/Hadamard/bell_state/length_circuits/l2qb.dat};
        \addplot+[brown,dashed,mark=*] table[ x={x},y={q4},col sep=comma, mark=* ] {data/simulator/1D/Hadamard/bell_state/length_circuits/l2qb.dat};
        \addplot+[black,dashed,mark=*] table[ x={x},y={q5},col sep=comma, mark=* ] {data/simulator/1D/Hadamard/bell_state/length_circuits/l2qb.dat};
\legend{SIM 2 qubits,
SIM 3 qubits,
SIM 4 qubits,
SIM 5 qubits,
hadamard test 2 qubits,
hadamard test 3 qubits,
hadamard test 4 qubits,
hadamard test 5 qubits
}
\end{axis}
\end{tikzpicture}
}
    \caption{circuit length }
    \label{fig:circuit_length}
\end{figure}
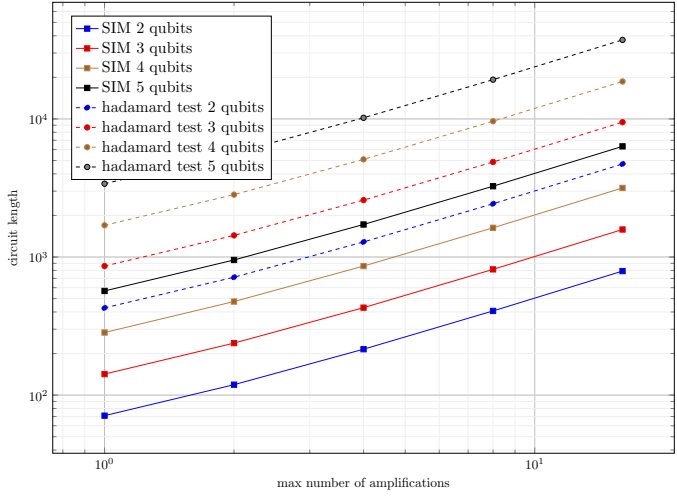
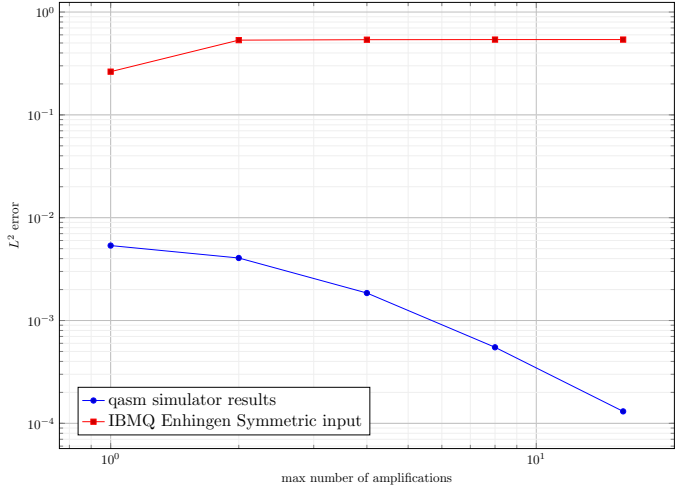
\begin{figure}
    \centering
    \resizebox{0.5\textwidth}{!}{
\begin{tikzpicture}
\begin{axis}[
    legend style={nodes={scale=1.3, transform shape}},
    xlabel=max number of amplifications,
    ylabel=$L^2$ error,
    xmode=log,
    ymode=log,
    grid = both,
    minor tick num = 1,
    major grid style = {lightgray},
    minor grid style = {lightgray!25},
    width = \textwidth,
    height = 0.75\textwidth,
    legend cell align = {left},
    legend pos = south west,
    ]
        \addplot table[ x={x},y expr=\thisrow{L2}^(1/2),col sep=comma, mark=* ] {\opbLt};
        \addplot table[ x={x},y expr=\thisrow{L2}^(1/2),col sep=comma, mark=* ] {\YopbLt};
\legend{qasm simulator results,
IBMQ Enhingen Symmetric input ,
IBMQ Enhingen Hadamard test
}
\end{axis}
\end{tikzpicture}
}
    \caption{real devices 2 qbits}
    \label{fig:my_label}
\end{figure}
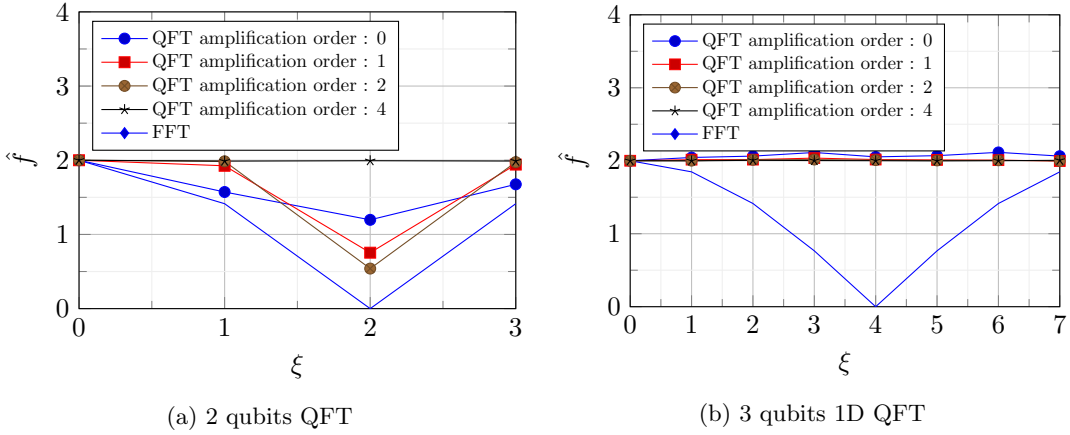
\begin{figure*}
    \centering
\begin{subfigure}{0.4\textwidth}
    \centering
    \resizebox{1\textwidth}{!}{
    \begin{tikzpicture}
\begin{axis}[
    legend style={nodes={scale=0.7, transform shape}},
    xlabel=$\xi$,
    ylabel=$\hat{f}$,
    xmin = 0, xmax = 3,
    ymin = 0, ymax = 4,
    xtick distance = 1,
    ytick distance = 1,
    grid = both,
    minor tick num = 1,
    major grid style = {lightgray},
    minor grid style = {lightgray!25},
    width = \textwidth,
    height = 0.75\textwidth,
    legend cell align = {left},
    legend pos = north west
]
 
        \addplot table[x={x},y expr=4*\thisrow{q0},col sep=comma, mark= *] {\Yopbt};
        \addplot table[x={x},y expr=4*\thisrow{q1},col sep=comma, mark= *] {\Yopbt};
        \addplot table[x={x},y expr=4*\thisrow{q2},col sep=comma, mark= *] {\Yopbt};
        \addplot table[x={x},y expr=4*\thisrow{q3},col sep=comma, mark= *] {\Yopbt};
        \addplot table[color=red, x={x},y expr=4*\thisrow{rv},col sep=comma, mark= ] {\Yopbt};
 
\legend{
    QFT amplification order : 0, 
    QFT amplification order : 1, 
    QFT amplification order : 2, 
    QFT amplification order : 4, 
    FFT
}
 
\end{axis}
 
\end{tikzpicture}}
    \caption{2 qubits QFT}
    \label{fig:real_SYM_2q}
\end{subfigure}
\begin{subfigure}{0.4\textwidth}
    \centering
    \resizebox{1\textwidth}{!}{
    \begin{tikzpicture}
\begin{axis}[
    legend style={nodes={scale=0.7, transform shape}},
    xlabel=$\xi$,
    ylabel=$\hat{f}$,
    xmin = 0, xmax = 7,
    ymin = 0, ymax =4,
    xtick distance = 1,
    ytick distance = 1,
    grid = both,
    minor tick num = 1,
    major grid style = {lightgray},
    minor grid style = {lightgray!25},
    width = \textwidth,
    height = 0.75\textwidth,
    legend cell align = {left},
    legend pos = north west
]
 
        \addplot table[x={x},y expr=8*\thisrow{q0},col sep=comma, mark= *] {\Yopbth};
        \addplot table[x={x},y expr=8*\thisrow{q1},col sep=comma, mark= *] {\Yopbth};
        \addplot table[x={x},y expr=8*\thisrow{q2},col sep=comma, mark= *] {\Yopbth};
        \addplot table[x={x},y expr=8*\thisrow{q3},col sep=comma, mark= *] {\Yopbth};
        \addplot table[color=red, x={x},y expr=8*\thisrow{rv},col sep=comma, mark= ] {\Yopbth};
 
\legend{
    QFT amplification order : 0, 
    QFT amplification order : 1, 
    QFT amplification order : 2, 
    QFT amplification order : 4, 
    FFT
}
 
\end{axis}
 
\end{tikzpicture}
}
    \caption{3 qubits 1D QFT}
    \label{fig:real_SYM_3q}
\end{subfigure}
    \caption{Bell state Input Quantum Fourier transform using sign determination methods on real device}
    \label{fig:real_SYM}
\end{figure*}
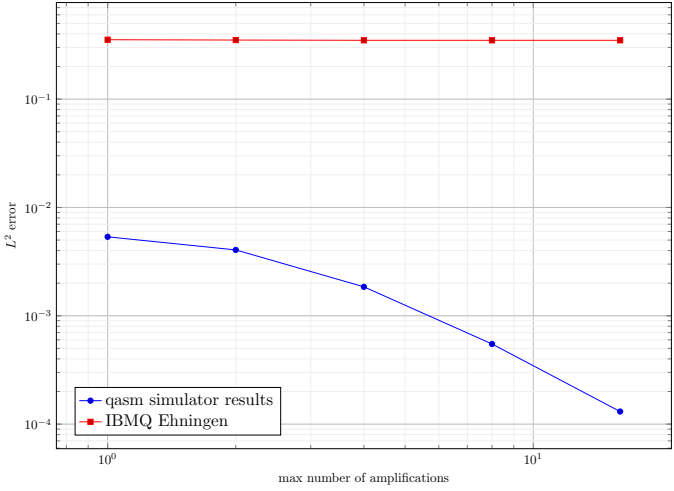
\begin{figure}
    \centering
    \resizebox{0.5\textwidth}{!}{
\begin{tikzpicture}
\begin{axis}[
    legend style={nodes={scale=1.3, transform shape}},
    xlabel=max number of amplifications,
    ylabel=$L^2$ error,
    xmode=log,
    ymode=log,
    grid = both,
    minor tick num = 1,
    major grid style = {lightgray},
    minor grid style = {lightgray!25},
    width = \textwidth,
    height = 0.75\textwidth,
    legend cell align = {left},
    legend pos = south west,
    ]
        \addplot table[ x={x},y expr=\thisrow{L2}^(1/2),col sep=comma, mark=* ] {\opbLt};
        \addplot table[ x={x},y expr=\thisrow{L2}^(1/2),col sep=comma, mark=* ] {\YopbLth};
\legend{qasm simulator results,
IBMQ Ehningen,
}
\end{axis}
\end{tikzpicture}
}
    \caption{sign determination real device 3 qbits}
    \label{fig:my_label}
\end{figure}
\begin{figure}
    \centering
    \resizebox{0.5\textwidth}{!}{
    \begin{tikzpicture}
\begin{axis}[
    legend style={nodes={scale=1.3, transform shape}},
    xlabel=$\xi$,
    ylabel=$\hat{f}$,
    xmin = 0, xmax = 3,
    ymin = 0, ymax = 4,
    xtick distance = 1,
    ytick distance = 1,
    grid = both,
    minor tick num = 1,
    major grid style = {lightgray},
    minor grid style = {lightgray!25},
    width = \textwidth,
    height = 0.75\textwidth,
    legend cell align = {left},
    legend pos = north west
]
 
        \addplot table[x={x},y expr=4*\thisrow{q0},col sep=comma, mark= *] {\YoHbt};
        \addplot table[x={x},y expr=4*\thisrow{q1},col sep=comma, mark= *] {\YoHbt};
        \addplot table[x={x},y expr=4*\thisrow{q2},col sep=comma, mark= *] {\YoHbt};
        \addplot table[color=red, x={x},y expr=4*\thisrow{rv},col sep=comma, mark= ] {\opbt};

\legend{
    QFT amplification order : 0, 
    QFT amplification order : 1, 
    QFT amplification order : 2, 
    FFT
}
 
\end{axis}
 
\end{tikzpicture}
}
    \caption{Bell state Input Quantum Fourier transform using Hadamard test on real devices (2 qubits)}
    \label{fig:my_label}
\end{figure}
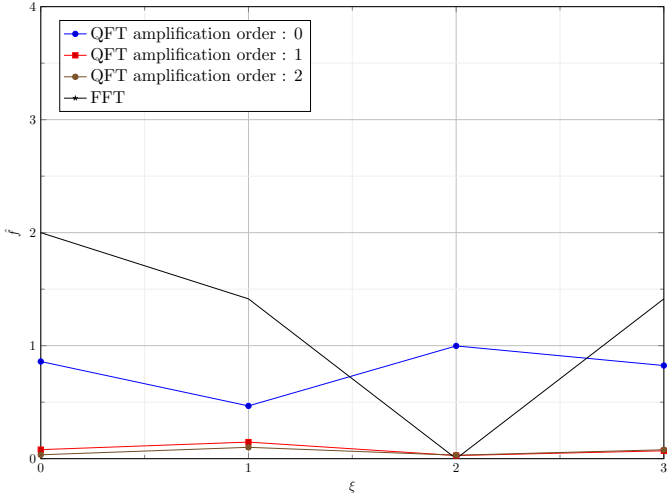
As a second step, the QFT in two dimensions with more realistic inputs have been implemented. Here the focus will only be on the sign determination algorithm. To test it here, either a random input vector (Figure \ref{fig:mirrored_random}) or a sinus input (Figure \ref{fig:mirrored_sinus}) is taken. Then a mirrored vector $f \in [0,1]^{2^{n+1}}$ is defined such as in Figure \ref{fig:mirrored_random} or in Figure \ref{fig:mirrored_sinus}. The main advantages of this input is that as it's symmetric which leads to exact Fourier coefficients estimation. Nevertheless, as it's a completely arbitrary input, the algorithm has to be linked with a an arbitrary state initialization algorithm which will increase the complexity of the quantum algorithm, but also the classical part will be more costly as all the gates for the initialization have to be calculated before. For this reason the algorithm will be only launched on simulator. As it can be seen on Figure \ref{fig:qft2d3}  and on Figure \ref{fig:qft2d3_sinus} the comportment observed in dimension 1 is extended in dimension 2, the local error decrease with the number of amplification. Nevertheless, the same problem remains that the larger the domain is, the bigger the amplification has to be to estimate usable Fourier coefficients. This problem may lead when it comes to the Strain calculation algorithm to error propagation which may disturb the algorithm convergence.
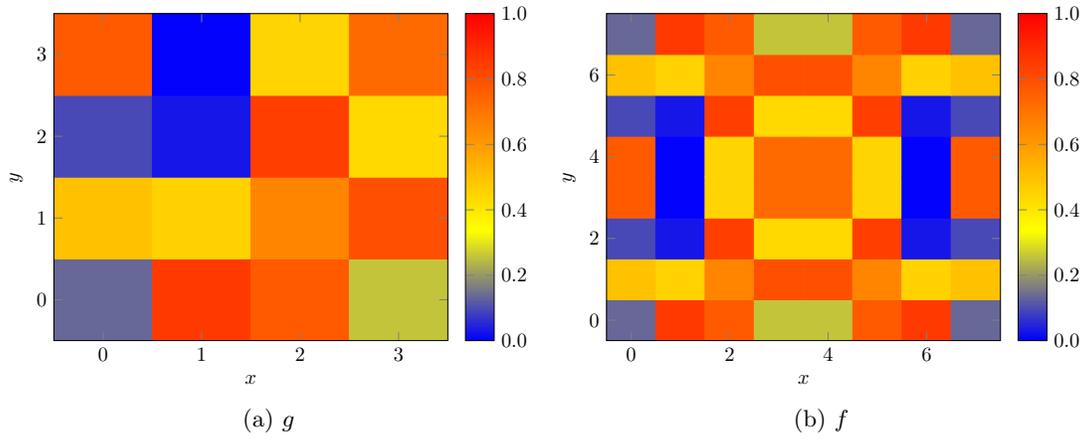
\begin{figure*}
    \centering
\begin{subfigure}{0.4\textwidth}
    \centering
        \resizebox{1\textwidth}{!}{
        \begin{tikzpicture}
          \begin{axis}[
            view={0}{90},   
            xlabel=$x$,
            ylabel=$y$,
            ymin=-0.5,ymax=3.5,
            xmin=-0.5,xmax=3.5,
            colorbar,
            colorbar style={
                yticklabel style={
                    /pgf/number format/.cd,
                    fixed,
                    precision=1,
                    fixed zerofill,
                },
            },
            enlargelimits=false,
            axis on top,
            point meta min=0,
            point meta max=1,
        ]
            \addplot [matrix plot*,point meta=explicit] file {data/simulator/2D/random3/ini_vec.dat};
          \end{axis}
        \end{tikzpicture}
        }
        \caption{$g$}
        \label{fig:my_label}
\end{subfigure}
\begin{subfigure}{0.4\textwidth}
    \centering
        \resizebox{1\textwidth}{!}{
        \begin{tikzpicture}
          \begin{axis}[
            view={0}{90},   
            xlabel=$x$,
            ylabel=$y$,
            ymin=-0.5,ymax=7.5,
            xmin=-0.5,xmax=7.5,
            colorbar,
            colorbar style={
                yticklabel style={
                    /pgf/number format/.cd,
                    fixed,
                    precision=1,
                    fixed zerofill,
                },
            },
            enlargelimits=false,
            axis on top,
            point meta min=0,
            point meta max=1,
        ]
            \addplot [matrix plot*,point meta=explicit] file {data/simulator/2D/random3/mirrored_ini_vec.dat};
          \end{axis}
        \end{tikzpicture}
        }
        \caption{$f$}
        \label{fig:my_label}
\end{subfigure}
    \caption{Random Mirrored input}
    \label{fig:mirrored_domain}
\end{figure*}

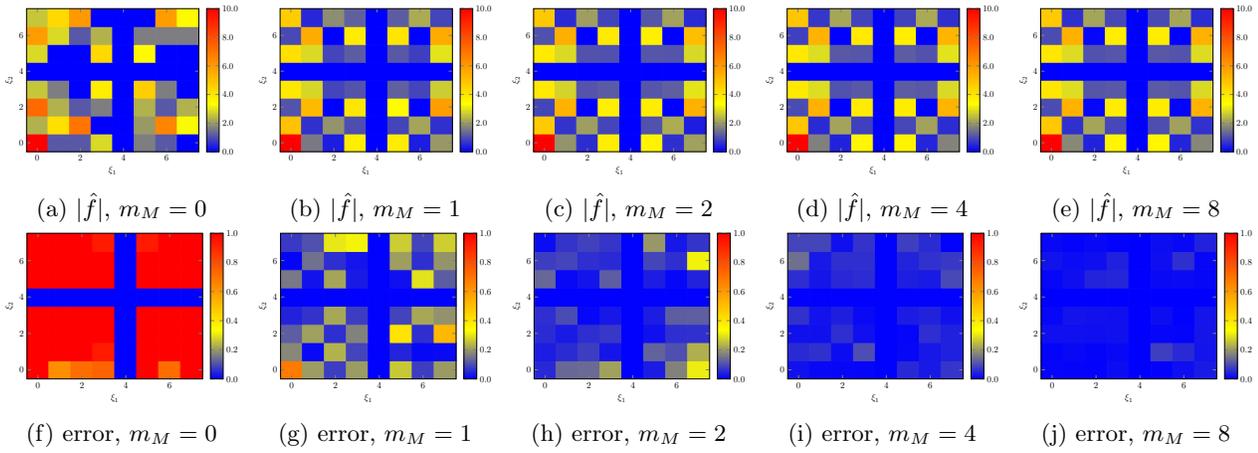
\begin{figure*}
    \begin{subfigure}{0.18\textwidth}
    \centering
    \resizebox{1\textwidth}{!}{
        \begin{tikzpicture}
          \begin{axis}[
            view={0}{90},   
            xlabel=$\xi_1$,
            ylabel=$\xi_2$,
            colorbar,
            colorbar style={
                yticklabel style={
                    /pgf/number format/.cd,
                    fixed,
                    precision=1,
                    fixed zerofill,
                },
            },
            enlargelimits=false,
            axis on top,
            point meta min=0,
            point meta max=10,
        ]
            \addplot [matrix plot*,point meta=explicit] file {data/simulator/2D/random3/qft/m_0.dat};
          \end{axis}
        \end{tikzpicture}
        }
        \caption{$|\hat{f}|$, $m_M=0$}
        \label{fig:mirrored_random}
    \end{subfigure}
    \begin{subfigure}{0.18\textwidth}
    \centering
    \resizebox{1\textwidth}{!}{
        \begin{tikzpicture}
          \begin{axis}[
            view={0}{90},   
            xlabel=$\xi_1$,
            ylabel=$\xi_2$,
            colorbar,
            colorbar style={
                yticklabel style={
                    /pgf/number format/.cd,
                    fixed,
                    precision=1,
                    fixed zerofill,
                },
            },
            enlargelimits=false,
            axis on top,
            point meta min=0,
            point meta max=10,
        ]
            \addplot [matrix plot*,point meta=explicit] file {data/simulator/2D/random3/qft/m_1.dat};
          \end{axis}
        \end{tikzpicture}
        }
        \caption{$|\hat{f}|$, $m_M=1$}
        \label{fig:my_label}
    \end{subfigure}
    \begin{subfigure}{0.18\textwidth}
    \centering
    \resizebox{1\textwidth}{!}{
        \begin{tikzpicture}
          \begin{axis}[
            view={0}{90},   
            xlabel=$\xi_1$,
            ylabel=$\xi_2$,
            colorbar,
            colorbar style={
                yticklabel style={
                    /pgf/number format/.cd,
                    fixed,
                    precision=1,
                    fixed zerofill,
                },
            },
            enlargelimits=false,
            axis on top,
            point meta min=0,
            point meta max=10,
        ]
            \addplot [matrix plot*,point meta=explicit] file {data/simulator/2D/random3/qft/m_2.dat};
          \end{axis}
        \end{tikzpicture}
        }
        \caption{$|\hat{f}|$, $m_M=2$}
        \label{fig:my_label}
    \end{subfigure}
    \begin{subfigure}{0.18\textwidth}
    \centering
    \resizebox{1\textwidth}{!}{
        \begin{tikzpicture}
          \begin{axis}[
            view={0}{90},   
            xlabel=$\xi_1$,
            ylabel=$\xi_2$,
            colorbar,
            colorbar style={
                yticklabel style={
                    /pgf/number format/.cd,
                    fixed,
                    precision=1,
                    fixed zerofill,
                },
            },
            enlargelimits=false,
            axis on top,
            point meta min=0,
            point meta max=10,
        ]
            \addplot [matrix plot*,point meta=explicit] file {data/simulator/2D/random3/qft/m_4.dat};
          \end{axis}
        \end{tikzpicture}
        }
        \caption{$|\hat{f}|$, $m_M=4$}
        \label{fig:my_label}
    \end{subfigure}
    \begin{subfigure}{0.18\textwidth}
    \centering
    \resizebox{1\textwidth}{!}{
        \begin{tikzpicture}
          \begin{axis}[
            view={0}{90},   
            xlabel=$\xi_1$,
            ylabel=$\xi_2$,
            colorbar,
            colorbar style={
                yticklabel style={
                    /pgf/number format/.cd,
                    fixed,
                    precision=1,
                    fixed zerofill,
                },
            },
            enlargelimits=false,
            axis on top,
            point meta min=0,
            point meta max=10,
        ]
            \addplot [matrix plot*,point meta=explicit] file {data/simulator/2D/random3/qft/m_8.dat};
          \end{axis}
        \end{tikzpicture}
        }
        \caption{$|\hat{f}|$, $m_M=8$}
        \label{fig:my_label}
    \end{subfigure}
    
    \begin{subfigure}{0.18\textwidth}
    \centering
    \resizebox{1\textwidth}{!}{
        \begin{tikzpicture}
          \begin{axis}[
            view={0}{90},   
            xlabel=$\xi_1$,
            ylabel=$\xi_2$,
            colorbar,
            colorbar style={
                yticklabel style={
                    /pgf/number format/.cd,
                    fixed,
                    precision=1,
                    fixed zerofill,
                },
            },
            enlargelimits=false,
            axis on top,
            point meta min=0,
            point meta max=1,
        ]
            \addplot [matrix plot*,point meta=explicit] file {data/simulator/2D/random3/L2/m_0.dat};
          \end{axis}
        \end{tikzpicture}
        }
        \caption{error, $m_M=0$}
        \label{fig:my_label}
    \end{subfigure}
    \begin{subfigure}{0.18\textwidth}
    \centering
        \resizebox{1\textwidth}{!}{
        \begin{tikzpicture}
          \begin{axis}[
            view={0}{90},   
            xlabel=$\xi_1$,
            ylabel=$\xi_2$,
            colorbar,
            colorbar style={
                yticklabel style={
                    /pgf/number format/.cd,
                    fixed,
                    precision=1,
                    fixed zerofill,
                },
            },
            enlargelimits=false,
            axis on top,
            point meta min=0,
            point meta max=1,
        ]
            \addplot [matrix plot*,point meta=explicit] file {data/simulator/2D/random3/L2/m_1.dat};
          \end{axis}
        \end{tikzpicture}
        }
        \caption{error, $m_M=1$}
        \label{fig:my_label}
    \end{subfigure}
    \begin{subfigure}{0.18\textwidth}
    \centering
        \resizebox{1\textwidth}{!}{
        \begin{tikzpicture}
          \begin{axis}[
            view={0}{90},   
            xlabel=$\xi_1$,
            ylabel=$\xi_2$,
            colorbar,
            colorbar style={
                yticklabel style={
                    /pgf/number format/.cd,
                    fixed,
                    precision=1,
                    fixed zerofill,
                },
            },
            enlargelimits=false,
            axis on top,
            point meta min=0,
            point meta max=1,
        ]
            \addplot [matrix plot*,point meta=explicit] file {data/simulator/2D/random3/L2/m_2.dat};
          \end{axis}
        \end{tikzpicture}
        }
        \caption{error, $m_M=2$}
        \label{fig:my_label}
    \end{subfigure}
    \begin{subfigure}{0.18\textwidth}
    \centering
        \resizebox{1\textwidth}{!}{
        \begin{tikzpicture}
          \begin{axis}[
            view={0}{90},   
            xlabel=$\xi_1$,
            ylabel=$\xi_2$,
            colorbar,
            colorbar style={
                yticklabel style={
                    /pgf/number format/.cd,
                    fixed,
                    precision=1,
                    fixed zerofill,
                },
            },
            enlargelimits=false,
            axis on top,
            point meta min=0,
            point meta max=1,
        ]
            \addplot [matrix plot*,point meta=explicit] file {data/simulator/2D/random3/L2/m_4.dat};
          \end{axis}
        \end{tikzpicture}
        }
        \caption{error, $m_M=4$}
        \label{fig:my_label}
    \end{subfigure}
    \begin{subfigure}{0.18\textwidth}
    \centering
        \resizebox{1\textwidth}{!}{
        \begin{tikzpicture}
          \begin{axis}[
            view={0}{90},   
            xlabel=$\xi_1$,
            ylabel=$\xi_2$,
            colorbar,
            colorbar style={
                yticklabel style={
                    /pgf/number format/.cd,
                    fixed,
                    precision=1,
                    fixed zerofill,
                },
            },
            enlargelimits=false,
            axis on top,
            point meta min=0,
            point meta max=1,
        ]
            \addplot [matrix plot*,point meta=explicit] file {data/simulator/2D/random3/L2/m_8.dat};
          \end{axis}
        \end{tikzpicture}
        }
        \caption{error, $m_M=8$}
        \label{fig:my_label}
    \end{subfigure}
    \caption{2D results on 3 qubits for random mirrored input}
    \label{fig:qft2d3}
\end{figure*}

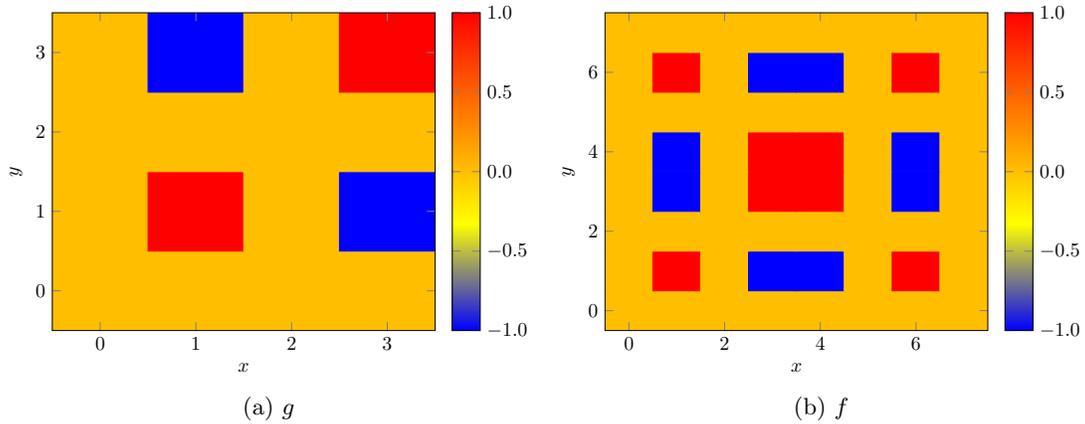
\begin{figure*}
    \centering
\begin{subfigure}{0.4\textwidth}
    \centering
        \resizebox{1\textwidth}{!}{
        \begin{tikzpicture}
          \begin{axis}[
            view={0}{90},   
            xlabel=$x$,
            ylabel=$y$,
            ymin=-0.5,ymax=3.5,
            xmin=-0.5,xmax=3.5,
            colorbar,
            colorbar style={
                yticklabel style={
                    /pgf/number format/.cd,
                    fixed,
                    precision=1,
                    fixed zerofill,
                },
            },
            enlargelimits=false,
            axis on top,
            point meta min=-1,
            point meta max=1,
        ]
            \addplot [matrix plot*,point meta=explicit] file {data/simulator/2D/sin3/ini_vec.dat};
          \end{axis}
        \end{tikzpicture}
        }
        \caption{$g$}
        \label{fig:my_label}
\end{subfigure}
\begin{subfigure}{0.4\textwidth}
    \centering
        \resizebox{1\textwidth}{!}{
        \begin{tikzpicture}
          \begin{axis}[
            view={0}{90},   
            xlabel=$x$,
            ylabel=$y$,
            ymin=-0.5,ymax=7.5,
            xmin=-0.5,xmax=7.5,
            colorbar,
            colorbar style={
                yticklabel style={
                    /pgf/number format/.cd,
                    fixed,
                    precision=1,
                    fixed zerofill,
                },
            },
            enlargelimits=false,
            axis on top,
            point meta min=-1,
            point meta max=1,
        ]
            \addplot [matrix plot*,point meta=explicit] file {data/simulator/2D/sin3/mirrored_ini_vec.dat};
          \end{axis}
        \end{tikzpicture}
        }
        \caption{$f$}
        \label{fig:my_label}
\end{subfigure}
    \caption{Mirrored sinus input}
    \label{fig:mirrored_sinus}
\end{figure*}

\begin{figure*}
    \begin{subfigure}{0.18\textwidth}
    \centering
    \resizebox{1\textwidth}{!}{
        \begin{tikzpicture}
          \begin{axis}[
            view={0}{90},   
            xlabel=$\xi_1$,
            ylabel=$\xi_2$,
            colorbar,
            colorbar style={
                yticklabel style={
                    /pgf/number format/.cd,
                    fixed,
                    precision=1,
                    fixed zerofill,
                },
            },
            enlargelimits=false,
            axis on top,
            point meta min=0,
            point meta max=10,
        ]
            \addplot [matrix plot*,point meta=explicit] file {data/simulator/2D/sin3/qft/m_0.dat};
          \end{axis}
        \end{tikzpicture}
        }
        \caption{$|\hat{f}|$, $m_M=0$}
        \label{fig:my_label}
    \end{subfigure}
    \begin{subfigure}{0.18\textwidth}
    \centering
    \resizebox{1\textwidth}{!}{
        \begin{tikzpicture}
          \begin{axis}[
            view={0}{90},   
            xlabel=$\xi_1$,
            ylabel=$\xi_2$,
            colorbar,
            colorbar style={
                yticklabel style={
                    /pgf/number format/.cd,
                    fixed,
                    precision=1,
                    fixed zerofill,
                },
            },
            enlargelimits=false,
            axis on top,
            point meta min=0,
            point meta max=10,
        ]
            \addplot [matrix plot*,point meta=explicit] file {data/simulator/2D/sin3/qft/m_1.dat};
          \end{axis}
        \end{tikzpicture}
        }
        \caption{$|\hat{f}|$, $m_M=1$}
        \label{fig:my_label}
    \end{subfigure}
    \begin{subfigure}{0.18\textwidth}
    \centering
    \resizebox{1\textwidth}{!}{
        \begin{tikzpicture}
          \begin{axis}[
            view={0}{90},   
            xlabel=$\xi_1$,
            ylabel=$\xi_2$,
            colorbar,
            colorbar style={
                yticklabel style={
                    /pgf/number format/.cd,
                    fixed,
                    precision=1,
                    fixed zerofill,
                },
            },
            enlargelimits=false,
            axis on top,
            point meta min=0,
            point meta max=10,
        ]
            \addplot [matrix plot*,point meta=explicit] file {data/simulator/2D/sin3/qft/m_2.dat};
          \end{axis}
        \end{tikzpicture}
        }
        \caption{$|\hat{f}|$, $m_M=2$}
        \label{fig:my_label}
    \end{subfigure}
    \begin{subfigure}{0.18\textwidth}
    \centering
    \resizebox{1\textwidth}{!}{
        \begin{tikzpicture}
          \begin{axis}[
            view={0}{90},   
            xlabel=$\xi_1$,
            ylabel=$\xi_2$,
            colorbar,
            colorbar style={
                yticklabel style={
                    /pgf/number format/.cd,
                    fixed,
                    precision=1,
                    fixed zerofill,
                },
            },
            enlargelimits=false,
            axis on top,
            point meta min=0,
            point meta max=10,
        ]
            \addplot [matrix plot*,point meta=explicit] file {data/simulator/2D/sin3/qft/m_4.dat};
          \end{axis}
        \end{tikzpicture}
        }
        \caption{$|\hat{f}|$, $m_M=4$}
        \label{fig:my_label}
    \end{subfigure}
    \begin{subfigure}{0.18\textwidth}
    \centering
    \resizebox{1\textwidth}{!}{
        \begin{tikzpicture}
          \begin{axis}[
            view={0}{90},   
            xlabel=$\xi_1$,
            ylabel=$\xi_2$,
            colorbar,
            colorbar style={
                yticklabel style={
                    /pgf/number format/.cd,
                    fixed,
                    precision=1,
                    fixed zerofill,
                },
            },
            enlargelimits=false,
            axis on top,
            point meta min=0,
            point meta max=10,
        ]
            \addplot [matrix plot*,point meta=explicit] file {data/simulator/2D/sin3/qft/m_8.dat};
          \end{axis}
        \end{tikzpicture}
        }
        \caption{$|\hat{f}|$, $m_M=8$}
        \label{fig:my_label}
    \end{subfigure}
    
    \begin{subfigure}{0.18\textwidth}
    \centering
    \resizebox{1\textwidth}{!}{
        \begin{tikzpicture}
          \begin{axis}[
            view={0}{90},   
            xlabel=$\xi_1$,
            ylabel=$\xi_2$,
            colorbar,
            colorbar style={
                yticklabel style={
                    /pgf/number format/.cd,
                    fixed,
                    precision=1,
                    fixed zerofill,
                },
            },
            enlargelimits=false,
            axis on top,
            point meta min=0,
            point meta max=1,
        ]
            \addplot [matrix plot*,point meta=explicit] file {data/simulator/2D/sin3/L2/m_0.dat};
          \end{axis}
        \end{tikzpicture}
        }
        \caption{error, $m_M=0$}
        \label{fig:my_label}
    \end{subfigure}
    \begin{subfigure}{0.18\textwidth}
    \centering
        \resizebox{1\textwidth}{!}{
        \begin{tikzpicture}
          \begin{axis}[
            view={0}{90},   
            xlabel=$\xi_1$,
            ylabel=$\xi_2$,
            colorbar,
            colorbar style={
                yticklabel style={
                    /pgf/number format/.cd,
                    fixed,
                    precision=1,
                    fixed zerofill,
                },
            },
            enlargelimits=false,
            axis on top,
            point meta min=0,
            point meta max=1,
        ]
            \addplot [matrix plot*,point meta=explicit] file {data/simulator/2D/sin3/L2/m_1.dat};
          \end{axis}
        \end{tikzpicture}
        }
        \caption{error, $m_M=1$}
        \label{fig:my_label}
    \end{subfigure}
    \begin{subfigure}{0.18\textwidth}
    \centering
        \resizebox{1\textwidth}{!}{
        \begin{tikzpicture}
          \begin{axis}[
            view={0}{90},   
            xlabel=$\xi_1$,
            ylabel=$\xi_2$,
            colorbar,
            colorbar style={
                yticklabel style={
                    /pgf/number format/.cd,
                    fixed,
                    precision=1,
                    fixed zerofill,
                },
            },
            enlargelimits=false,
            axis on top,
            point meta min=0,
            point meta max=1,
        ]
            \addplot [matrix plot*,point meta=explicit] file {data/simulator/2D/sin3/L2/m_2.dat};
          \end{axis}
        \end{tikzpicture}
        }
        \caption{error, $m_M=2$}
        \label{fig:my_label}
    \end{subfigure}
    \begin{subfigure}{0.18\textwidth}
    \centering
        \resizebox{1\textwidth}{!}{
        \begin{tikzpicture}
          \begin{axis}[
            view={0}{90},   
            xlabel=$\xi_1$,
            ylabel=$\xi_2$,
            colorbar,
            colorbar style={
                yticklabel style={
                    /pgf/number format/.cd,
                    fixed,
                    precision=1,
                    fixed zerofill,
                },
            },
            enlargelimits=false,
            axis on top,
            point meta min=0,
            point meta max=1,
        ]
            \addplot [matrix plot*,point meta=explicit] file {data/simulator/2D/sin3/L2/m_4.dat};
          \end{axis}
        \end{tikzpicture}
        }
        \caption{error, $m_M=4$}
        \label{fig:my_label}
    \end{subfigure}
    \begin{subfigure}{0.18\textwidth}
    \centering
        \resizebox{1\textwidth}{!}{
        \begin{tikzpicture}
          \begin{axis}[
            view={0}{90},   
            xlabel=$\xi_1$,
            ylabel=$\xi_2$,
            colorbar,
            colorbar style={
                yticklabel style={
                    /pgf/number format/.cd,
                    fixed,
                    precision=1,
                    fixed zerofill,
                },
            },
            enlargelimits=false,
            axis on top,
            point meta min=0,
            point meta max=1,
        ]
            \addplot [matrix plot*,point meta=explicit] file {data/simulator/2D/sin3/L2/m_8.dat};
          \end{axis}
        \end{tikzpicture}
        }
        \caption{error, $m_M=8$}
        \label{fig:my_label}
    \end{subfigure}
    \caption{2D results on 3 qubits for sinus mirrored input}
    \label{fig:qft2d3_sinus}
\end{figure*}
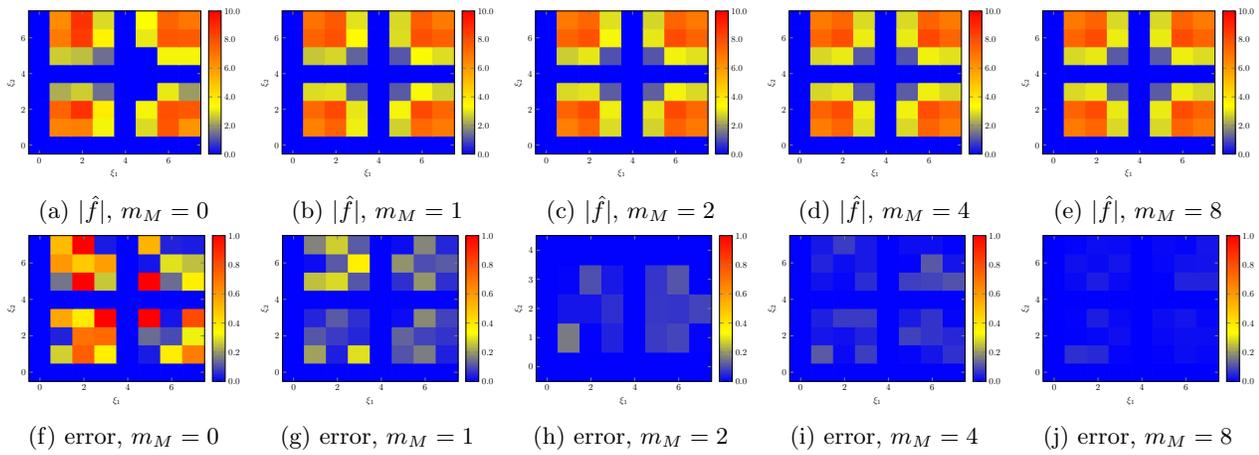

\section{QFT-based Homogenization}

\subsection{QFT based Homogenization algorithm}

The aim now is to integrate this different algorithms in an Homogenization code. To do so, data encoding has to be reworked. As the different Fourier transforms are performed on $\boldsymbol{\tau}$ which in 2D has 3 components. The quantum state encoded in the register has to be  \begin{equation}
    \ket{\boldsymbol{\tau}}{}_{n_1,n_2,n_d}=\sum_{l=1}^{d} \sum_{k=0}^{N_1-1}\sum_{j=0}^{N_2-1}\tau_{lkj}\ket{k}{}_{n_1} \ket{j}{}_{n_2}\ket{l}_{n_d},
\end{equation} where $\ket{\boldsymbol{\tau}}_{n_1,n_2,n_d}$ stands for $\frac{\boldsymbol{\tau}}{\|\boldsymbol{\tau}\|}$  with $N_1 = 2^{n_1}$ and $N_2 = 2^{n_2}$ where $n_1$ and $n_2$ are the number of qubits used to encode the data in the two dimensions. Here $d=3$ but the number of qubits used to implement the different components is $n_d=2$ using the fact that $\boldsymbol{\tau}_4 = \boldsymbol{0}$.

Then in the case of periodic Boundary conditions Algorithm \ref{alg:FFT_based_alg} can be rewritten 
\begin{algorithm}[H]
\caption{QFT based STRAIN}\label{alg:QFT_based_alg}
\begin{algorithmic}
\Require $\boldsymbol{E},\boldsymbol{C},\boldsymbol{M}$, 
\State $\boldsymbol{\varepsilon}^0(\boldsymbol{x})=\boldsymbol{E}, \forall \boldsymbol{x} \in Y,$
\State $k=1$
\While{convergence is not reached}
\State $\boldsymbol{\tau}^k = (\boldsymbol{C}(\boldsymbol{x})-\boldsymbol{C}^0):\boldsymbol{\varepsilon}^k(\boldsymbol{x}) $
\State $\boldsymbol{\hat{\tau}}^k(\boldsymbol{\xi})=\text{Hadamard QFT}({\boldsymbol{\tau}}^k(\boldsymbol{x}),\boldsymbol{M}),$
\State $\boldsymbol{\hat{\tau}}^{k+1}(\boldsymbol{\xi})=\boldsymbol{\hat{\Gamma}}^0(\boldsymbol{\xi}):\boldsymbol{\hat{\tau}}^k(\boldsymbol{\xi}),$
\State $\boldsymbol{\varepsilon}^{k+1}(\boldsymbol{x})=\text{Hadamard IQFT}(\boldsymbol{\hat{\tau}}^{k+1}(\boldsymbol{\xi}),\boldsymbol{M}),$
\State $\boldsymbol{\varepsilon}^{k+1}(\boldsymbol{x})=\boldsymbol{E}-\boldsymbol{\varepsilon}^{k+1}(\boldsymbol{x})$
\State $k=k+1$,
\State convergence test,
\EndWhile
\State $\boldsymbol{\sigma}(\boldsymbol{x})=\boldsymbol{C}(\boldsymbol{x}):\boldsymbol{\varepsilon}^{k}(\boldsymbol{x}),$
\return $\boldsymbol{\sigma}$
\end{algorithmic}
\end{algorithm}
Parameter m has to be chosen wisely in order to ensure general convergence. Differently said, MLQAE accuracy has to be big enough to ensure that global limiting error is coming from the gradient descent and not from the coefficients estimation.

In the case of symmetric input data, the algorithm is transformed similarly but using instead Algorithm \ref{alg:symmetric_QFT} \begin{algorithm}[H]
\caption{Symmetric QFT based STRAIN}\label{alg:Symmetric_QFT_based_alg}
\begin{algorithmic}
\Require $\boldsymbol{E},\boldsymbol{C},\boldsymbol{M}$, 
\State $\boldsymbol{\varepsilon}^0(\boldsymbol{x})=\boldsymbol{E}, \forall \boldsymbol{x} \in Y,$
\State $k=1$
\While{convergence is not reached}
\State $\boldsymbol{\tau}^k = (\boldsymbol{C}(\boldsymbol{x})-\boldsymbol{C}^0):\boldsymbol{\varepsilon}^k(\boldsymbol{x}) $
\State $\boldsymbol{\hat{\tau}}^k(\boldsymbol{\xi})=\text{Symmetric QFT}({\boldsymbol{\tau}}^k(\boldsymbol{x}),\boldsymbol{M}),$
\State $\boldsymbol{\hat{\tau}}^{k+1}(\boldsymbol{\xi})=\boldsymbol{\hat{\Gamma}}^0(\boldsymbol{\xi}):\boldsymbol{\hat{\tau}}^k(\boldsymbol{\xi}),$
\State $\boldsymbol{\varepsilon}^{k+1}(\boldsymbol{x})=\text{Symmetric IQFT}(\boldsymbol{\hat{\tau}}^{k+1}(\boldsymbol{\xi}),\boldsymbol{M}),$
\State $\boldsymbol{\varepsilon}^{k+1}(\boldsymbol{x})=\boldsymbol{E}-\boldsymbol{\varepsilon}^{k+1}(\boldsymbol{x})$
\State $k=k+1$,
\State convergence test,
\EndWhile
\State $\boldsymbol{\sigma}(\boldsymbol{x})=\boldsymbol{C}(\boldsymbol{x}):\boldsymbol{\varepsilon}^{k}(\boldsymbol{x}),$
\return $\boldsymbol{\sigma}$
\end{algorithmic}
\end{algorithm}

Finally a we are mainly interested into the effective overall material response $\boldsymbol{C^*}$ the Algorithm \ref{alg:QFT_based_sign} is transformed into the following algorithm (using de Voigt notations) \begin{algorithm}[H]
\caption{QFT-based Homogenization}
\begin{algorithmic}
\Require $\boldsymbol{C},\boldsymbol{M}$
\For{$j=1,3$}
\State $\boldsymbol{E}=\boldsymbol{0}$
\State $\boldsymbol{E}_{j}=\boldsymbol{1}$
\State $\boldsymbol{\sigma}= \text{QFT-based STRAIN}(\boldsymbol{E},\boldsymbol{C},\boldsymbol{M})$
\For{$k=1,3$}
\State $c^*_{jk}=\langle\boldsymbol{\sigma}_k\rangle_Y$
\EndFor
\EndFor
\return $\boldsymbol{c^*}$.
\end{algorithmic}
\end{algorithm}

\subsubsection{Results }

\begin{figure*}
\begin{subfigure}[b]{0.5\textwidth}
    \centering
    \includegraphics[width=0.89\textwidth]{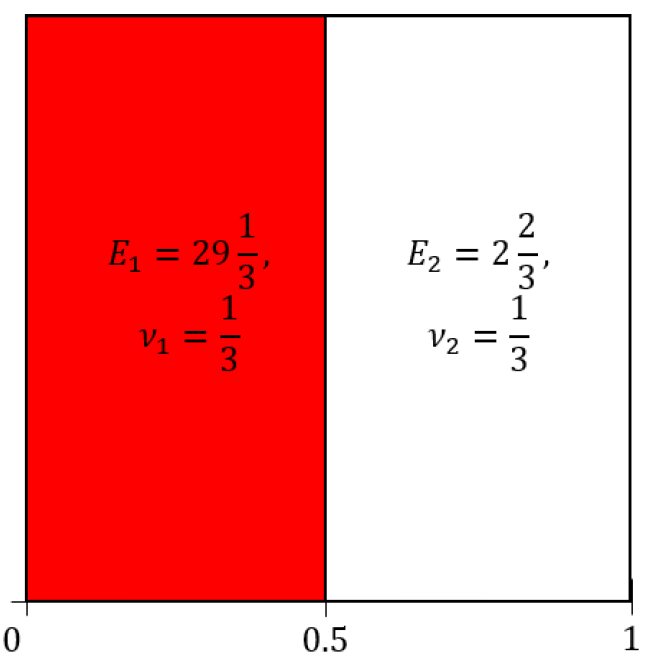}
    \caption{Laminate Geometry}
    \label{fig:laminate}
\end{subfigure}
\begin{subfigure}[b]{0.5\textwidth}
    \centering
    \includegraphics[width=1\textwidth]{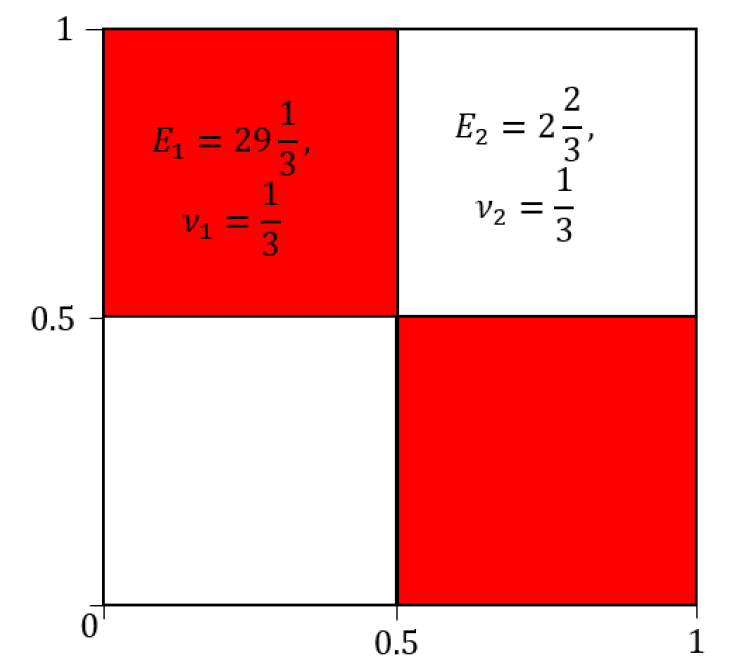}
    \caption{checkerboard Geometry}
    \label{fig:my_label}
\end{subfigure}
    \caption{Geometry}
    \label{fig:Geometries}
\end{figure*}

In this section we implemented and tested a QFT based homogenization solver for different geometries (Figure \ref{fig:Geometries}). The first one used is a laminate geometry which has for advantage to give under certain parameters the exact solution (given by \cite{https://doi.org/10.48550/arxiv.2204.13624}) after 1 iteration of the Strain algorithm. This advantage avoid the propagation of the QFT error for a sufficient number of amplification, and moreover allows to compare our solution with the exact one given with the FFT based algorithm. The second one is the checkerboard geometry which is a bit more complex and leads to an approximate solution with the algorithms, this will be important to estimate the influence of QFT error propagation on the algorithm convergence.
What is done here is to calculate the average strain field on a mirrored geometry to determine an approximation of the general stiffness matrix of the geometry under PMUBC. As stipulated for the laminate geometry the stiffness matrix is the exact one and looks like 
\begin{equation}\boldsymbol{c^*}_{\text{laminate}}=
    \begin{pmatrix}
    7.333333 & 3.666666 & 0 \\
     3.6666666 & 19.83333 & 0 \\
     0 & 0 & 1.8333333 \\
\end{pmatrix}.
\end{equation}
We tried to obtain a similar result with algorithm \ref{alg:QFT_based_alg} launched with a QASM simulator for quantum computations. The maximum number of amplification used is $2^{6}$ for a number of shots of 1000. The stiffness matrix obtained is
\begin{equation}\boldsymbol{C^*}_{\text{laminate}}=
    \begin{pmatrix}
    7.334 & 3.667 & 10e-11 \\
     3.667 & 19.832 & 10e-11 \\
     10e-11 & 10e-11 & 1.832 
\end{pmatrix}.
\end{equation}
It can be noticed that every component as a $1e-3$ accuracy but when it comes to null elements, the accuracy increase widely. This can be explained by the fact that a component equal to zero on a quantum computer has no chance to appear in the outcome measured which leads to an almost perfect result for this component only noised by the most likelihood method but still negligible. \\

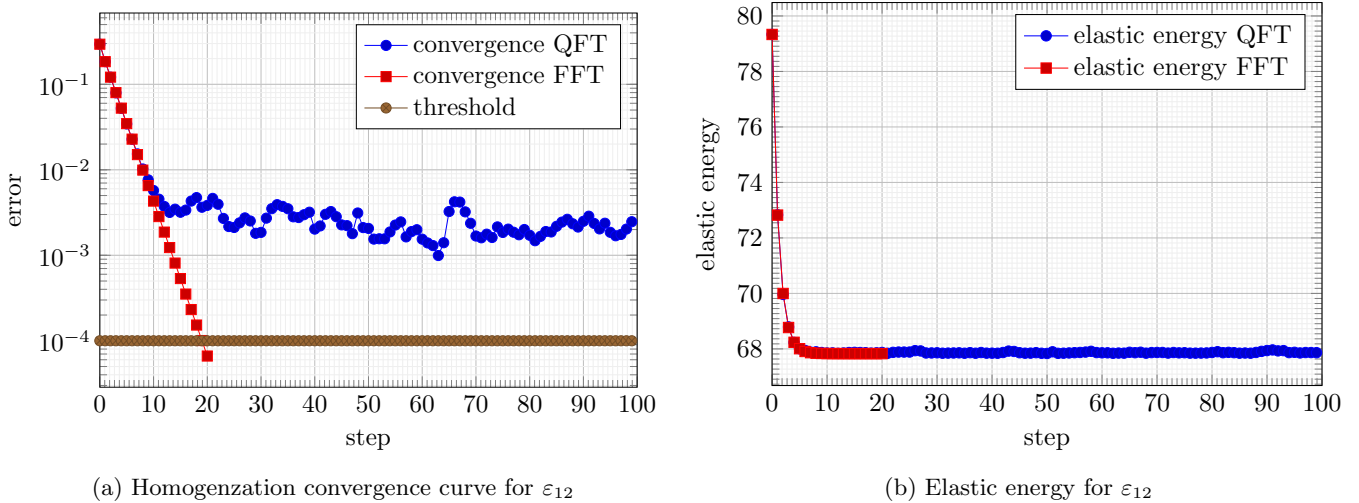
\begin{figure*}
\begin{subfigure}[b]{0.5\textwidth}
    \centering
    \resizebox{\textwidth}{!}{
\begin{tikzpicture}
\begin{axis}[
    legend style={nodes={scale=1, transform shape}},
    xlabel=step,
    ylabel=error,
    ymode=log,
    grid = both,
    xmin = 0, xmax = 100,
    minor tick num = 10,
    xtick distance = 10,
    major grid style = {lightgray},
    minor grid style = {lightgray!25},
    width = \textwidth,
    height = 0.75\textwidth,
    legend cell align = {left},
    legend pos = north east,
]
        \addplot table[ x={x},y={qft},col sep=comma, mark=* ] {data/simulator/2D/homogenization/qft_hom_errors.dat};
        \addplot table[ x={x},y={qft},col sep=comma, mark=* ] {data/simulator/2D/homogenization/fft_hom_errors.dat};
        \addplot table[ x={x},y={thr},col sep=comma, mark=* ] {data/simulator/2D/homogenization/qft_hom_errors.dat};

\legend{
convergence QFT,
convergence FFT,
threshold,
}
\end{axis}
\end{tikzpicture}
}
    \caption{Homogenzation convergence curve for $\varepsilon_{12}$}
    \label{fig:Hom_cv_error}
\end{subfigure}
\begin{subfigure}[b]{0.5\textwidth}
    \centering
        \resizebox{\textwidth}{!}{
\begin{tikzpicture}
\begin{axis}[
    legend style={nodes={scale=1, transform shape}},
    xlabel=step,
    ylabel=elastic energy,
    grid = both,
    xmin = 0, xmax = 100,
    minor tick num = 10,
    xtick distance = 10,
    major grid style = {lightgray},
    minor grid style = {lightgray!25},
    width = \textwidth,
    height = 0.75\textwidth,
    legend cell align = {left},
    legend pos = north east,
]
        \addplot table[ x={x},y={qft},col sep=comma, mark=* ] {data/simulator/2D/homogenization/qft_hom_energy.dat};
        \addplot table[ x={x},y={qft},col sep=comma, mark=* ] {data/simulator/2D/homogenization/fft_hom_energy.dat};

\legend{
elastic energy QFT,
elastic energy FFT,
}
\end{axis}
\end{tikzpicture}
}
    \caption{Elastic energy for $\varepsilon_{12}$}
    \label{fig:ee_homogeniation}
\end{subfigure}
    \caption{checkerboard homogenization with maximum amplification $m_M=2^6$}
    \label{fig:Sh_hom}
\end{figure*}

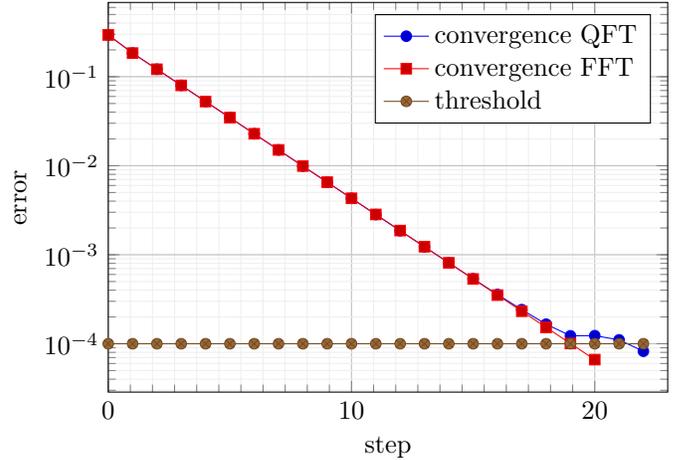
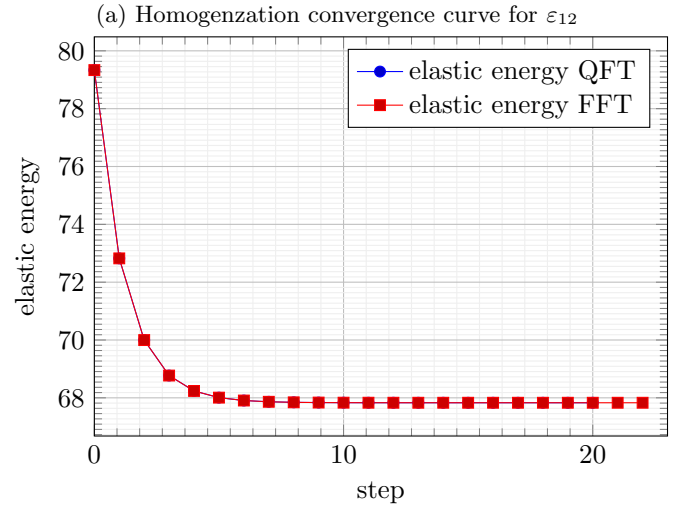
\begin{figure}
\begin{subfigure}[b]{0.5\textwidth}
    \centering
    \resizebox{\textwidth}{!}{
\begin{tikzpicture}
\begin{axis}[
    legend style={nodes={scale=1, transform shape}},
    xlabel=step,
    ylabel=error,
    ymode=log,
    grid = both,
    xmin = 0, xmax = 23,
    minor tick num = 10,
    xtick distance = 10,
    major grid style = {lightgray},
    minor grid style = {lightgray!25},
    width = \textwidth,
    height = 0.75\textwidth,
    legend cell align = {left},
    legend pos = north east,
]
        \addplot table[ x={x},y={qft},col sep=comma, mark=* ] {data/simulator/2D/homogenization/qft_hom_err_m9.dat};
        \addplot table[ x={x},y={qft},col sep=comma, mark=* ] {data/simulator/2D/homogenization/fft_hom_errors.dat};
        \addplot table[ x={x},y={th},col sep=comma, mark=* ] {data/simulator/2D/homogenization/qft_hom_err_m9.dat};

\legend{
convergence QFT,
convergence FFT,
threshold,
}
\end{axis}
\end{tikzpicture}
}
    \caption{Homogenzation convergence curve for $\varepsilon_{12}$}
    \label{fig:Hom_cv_error_m9}
\end{subfigure}
\begin{subfigure}[b]{0.5\textwidth}
    \centering
        \resizebox{\textwidth}{!}{
\begin{tikzpicture}
\begin{axis}[
    legend style={nodes={scale=1, transform shape}},
    xlabel=step,
    ylabel=elastic energy,
    grid = both,
    xmin = 0, xmax = 23,
    minor tick num = 10,
    xtick distance = 10,
    major grid style = {lightgray},
    minor grid style = {lightgray!25},
    width = \textwidth,
    height = 0.75\textwidth,
    legend cell align = {left},
    legend pos = north east,
]
        \addplot table[ x={x},y={qft},col sep=comma, mark=* ] {data/simulator/2D/homogenization/fft_hom_energy.dat};
        \addplot table[ x={x},y={qft},col sep=comma, mark=* ] {data/simulator/2D/homogenization/qft_hom_el_m9.dat};

\legend{
elastic energy QFT,
elastic energy FFT,
}
\end{axis}
\end{tikzpicture}
}
    \caption{Elastic energy for $\varepsilon_{12}$}
    \label{fig:ee_homogeniation_m9}
\end{subfigure}
    \caption{checkerboard homogenization with maximum amplification $m_M=2^9$}
    \label{fig:Sh_hom_m9}
\end{figure}

For the checkerboard geometry, the convergence is not reached in one iteration so a wise choice of convergence threshold has to be chosen. The literature tells us \cite{MOULINEC199869} that a convergence parameter of $1e-4$ will ensure an accurate enough solution to be physically acceptable which leads to the stiffness tensor
\begin{equation}\boldsymbol{C^*}_{\text{checkerboard}}=
    \begin{pmatrix}
    8.876 & 5.209 & 2.169 \\
     5.209 & 8.876 & 2.169 \\
     2.168 & 2.168 & 3.377 \\
\end{pmatrix}.
\end{equation}
Nevertheless, using this convergence treshold with the same MLQAE parameters is not sufficient as it doesn't lead to a convergence of the algorithm as it can be seen on figure \ref{fig:Hom_cv_error}. After 100 iterations the system is still not converged and gives the following 
\begin{equation}\boldsymbol{C^*}_{\text{checkerboard}}=
    \begin{pmatrix}
    8.849 & 5.188 & 2.169 \\
     5.170 & 8.869 & 2.171 \\
     2.180 & 2.173 & 3.377 \\
\end{pmatrix}.
\end{equation}
This homogenization stiffness matrix cannot really be used in practice as it's not physically acceptable due to this non-symmetric character with an error of $1e-2$. This difference can be explained by the fact that the MLQAE is adding noise at every step in the system which perturb the convergence. To hope for better results MLQAE number of amplifications should be increased in order to have noise negligible compared to the the convergence threshold. To do so, we would like to have MLQAE method accuracy reaching the classical computational machine error single precision ($1e-7$). Therefore,the maximum number of amplification is obtained using formula \ref{eqn:amplifcation_number} giving a maximum number of amplification around $2^{11}$. Nevertheless for the checkerboard geometry it seems that $2^9$ amplifications are enough to converge. The results for this experiment are shown in Figure \ref{fig:Sh_hom_m9}. It can be seen that the convergence is this time reached after 23 steps, this longer number of steps can be explained by the fact that noise will make a little deviation from solution as it sums up, and it then require few more steps to converge. The effective stiffness matrix is given by
\begin{equation}\boldsymbol{C^*}_{\text{checkerboard}}=
    \begin{pmatrix}
    8.877 & 5.210 & 2.168 \\
     5.209 & 8.877 & 2.168 \\
     2.168 & 2.168 & 3.376 \\
\end{pmatrix}.
\end{equation}
Similarly to the results obtained with the laminate geometry, the tensor is almost symmetric (with $1e-3$ error) which makes it usable for physical problems.  

\section{Conclusion}
In this paper we proposed an algorithm able to perform an hybrid Quantum Fourier Transform based Homogenization, we based this latter on the Moulinec and Suquet \cite{MOULINEC199869} gradient method which use a classical Fourier Transform. As a first step we implemented the Quantum Fourier Transform in different configurations (multi dimensional, applied to tensors) in order to have a functional QFT applicable to the Homogenization problem. Secondly we solved the problematic of the data transfer between classical computers and Quantum computers. The initialization part has been solved by taking one of the current most efficient algorithm made by Mottonen and al. \cite{mottonen2004transformation} allowing to have an initialization circuit growing exponentially with the number of qubits. The major problematic was coming from the reading out part and the fact that we had to estimate complex coefficients, in our researches we found two ways to do it, both based on the estimation of the amplitude of a state. Luckily methods for amplitude estimation has been widely developed in the last years especially by Suzuki and al. \cite{suzuki2020amplitude} using amplitude amplification and most likelihood methods to estimate the amplitude. As a little step forward on this subject we made an experimental analysis of the MLQAE method in order to have a way to determine the number of amplification necessary for a given accuracy. Our first way to estimate the Fourier coefficients was to use the MLQAE method mixed with an Hadamard test, allowing to estimate in two times the real part and the imaginary part of each coefficients. Nevertheless this method inherits the major problem of Hadamard test which is the presence of a complex controlled operator resulting in very deep circuits that are not suitable for the current noisy real quantum computers. In order to solve this problem we developed a second method based on the convenience of having a symmetric real input in the estimation of the Fourier coefficients. This method is based on the curious analogy between the resolution of the Limpman Schwinger equation with periodic boundary conditions on mirror domain and the resolution of the said equation with mixed boundary conditions studied first by Kabel  \cite{kabel2016mixed} and then developed recently by the latter and Grimm-strele \cite{GrimmStrele2021FFTB}. This would allow us to only deal with complex coefficients but of known phase, the resulting problem would then only be the determination of the sign of each coefficients. To do so we developed a method using the properties of complex exponential Fourier transform turning a sign determination problem into an amplitude comparison problem. This allowed us to have way shorter circuits but to the cost of having to consecutive MLQAE methods where one is depending on the other. Using this last method we then developed an Hybrid QFT based Homogenization solver that we tested on trivial geometries and that give similar results than the FFT based algorithm. Despite the apparent unsuitability of our solutions on real devices due to the noise, we hope for a future development of less noisy qubit using the advantage of noise correction methods or engineering tools such as IBM mapomatic to optimize circuits building in function of qubits noise. 

\section*{Acknowledgments} 

\emph{This work was supported by the project AnQuC-3 of the Competence Center Quantum Computing Rhineland-Palatinate
(Germany).}

\bibliographystyle{ieeetr} 
\bibliography{bibli.bib}

\end{document}